\documentclass[preprint]{elsarticle}  

\makeatletter
\def\ps@pprintTitle{%
  \let\@oddhead\@empty
  \let\@evenhead\@empty
  \def\@oddfoot{\centerline{\thepage}}%
  \let\@evenfoot\@oddfoot}
\makeatother

\usepackage{etoolbox}
\patchcmd{\MaketitleBox}%
  {\footnotesize\itshape\elsaddress\par\vskip36pt}%
  {\footnotesize\itshape\elsaddress\par\parbox[b][36pt]{\linewidth}{\vfill\hfill\textnormal{\today}\hfill\null\vfill}}%
  {}{}
\patchcmd{\pprintMaketitle}%
  {\footnotesize\itshape\elsaddress\par\vskip36pt}%
  {\footnotesize\itshape\elsaddress\par\parbox[b][36pt]{\linewidth}{\vfill\hfill\textnormal{\today}\hfill\null\vfill}}%
  {}{}

\usepackage{float}
\usepackage{amsfonts} 
\usepackage{caption}
\usepackage{siunitx}

\usepackage[
            colorlinks=true,%
            breaklinks=true,%
            linkcolor=blue,
            urlcolor=blue,%
            citecolor=blue,%
            pdftitle={Title of paper}, 
            pdfkeywords={Keywords},	
            pdfauthor={Authors},		
            bookmarksopen=false,
            pdfpagemode=UseNone]{hyperref}

\usepackage[margin = 1.2in]{geometry}
\usepackage[T1]{fontenc}
\usepackage[english]{babel}
\usepackage[utf8]{inputenc}
\usepackage{mathtools}
\usepackage{amsmath}
\usepackage{amsfonts}
\usepackage{amsthm}
\usepackage{amssymb}
\usepackage{array}
\usepackage{algorithm} 
\usepackage{algorithmicx}
\usepackage{algpseudocode}
\usepackage{subcaption}
\usepackage{siunitx}
\usepackage{arydshln}
\usepackage{nomencl}
\usepackage{framed}  
\usepackage{multicol}  

\usepackage{color}
\usepackage{xcolor}
\usepackage{url}
\usepackage{cleveref}
\usepackage{graphicx}
\usepackage{subcaption}
\usepackage{breakcites}
\usepackage{soul} 
\usepackage{enumitem}
\usepackage{tikz}
\usetikzlibrary{arrows.meta, positioning}
\usetikzlibrary{patterns}
\usepackage{pgfplots}
\pgfplotsset{compat=1.16}

{\noindent {\textbf{Proof}:} }%
{\hfill $\Box$ \\[1ex] }

\newcommand{\bit}{\begin{itemize}}
\newcommand{\eit}{\end{itemize}}
\newcommand{\ben}{\begin{enumerate}}
\newcommand{\een}{\end{enumerate}}

\newcommand {\real} {\mathbb{R}}
\newcommand {\nat} {\mathbb{N}}

\newcommand{\bA}{\ensuremath{\mathbf{A}}}
\newcommand{\bB}{\ensuremath{\mathbf{B}}}

\newcommand{\bD}{\ensuremath{\mathbf{D}}}

\newcommand{\bH}{\ensuremath{\mathbf{H}}}
\newcommand{\bI}{\ensuremath{\mathbf{I}}}

\newcommand{\bO}{\ensuremath{\mathbf{O}}}

\newcommand{\bS}{\ensuremath{\mathbf{S}}}

\newcommand{\bU}{\ensuremath{\mathbf{U}}}
\newcommand{\bV}{\ensuremath{\mathbf{V}}}
\newcommand{\bW}{\ensuremath{\mathbf{W}}}
\newcommand{\bX}{\ensuremath{\mathbf{X}}}

\newcommand{\ba}{\ensuremath{\mathbf{a}}}

\newcommand{\bc}{\ensuremath{\mathbf{c}}}

\renewcommand{\bf}{\ensuremath{\mathbf{f}}}

\newcommand{\bs}{\ensuremath{\mathbf{s}}}

\newcommand{\bu}{\ensuremath{\mathbf{u}}}

\newcommand{\cE}{\ensuremath{\mathcal{E}}}

\newcommand{\cM}{\ensuremath{\mathcal{M}}}

\begin{document}

\title{Parametric Operator Inference to Simulate the Purging Process in Semiconductor Manufacturing}

\author[1]{Seunghyon Kang}    
\author[2]{Hyeonghun Kim\corref{cor1}}     
\ead{hyk049@ucsd.edu}
\author[2]{Boris Kramer}
\address[1]{Memory Defect Science \& Engineering Team, Samsung Electronics Co., Ltd, Hwaseong-si, Republic of Korea} 
\address[2]{Department of Mechanical and Aerospace Engineering, University of California San Diego, CA, United States}

\cortext[cor1]{Corresponding author}

\begin{abstract}
This work presents the application of parametric Operator Inference (OpInf)---a nonintrusive reduced-order modeling (ROM) technique that learns a low-dimensional representation of a high-fidelity model---to the numerical model of the purging process in semiconductor manufacturing. Leveraging the data-driven nature of the OpInf framework, we aim to forecast the flow field within a plasma-enhanced chemical vapor deposition~(PECVD) chamber using computational fluid dynamics (CFD) simulation data. Our model simplifies the system by excluding plasma dynamics and chemical reactions, while still capturing the key features of the purging flow behavior.
The parametric OpInf framework learns nine ROMs based on varying argon mass flow rates at the inlet and different outlet pressures. It then interpolates these ROMs to predict the system's behavior for 25 parameter combinations, including 16 scenarios that are not seen in training.
The parametric OpInf ROMs, trained on 36\% of the data and tested on 64\%, demonstrate accuracy across the entire parameter domain, with a maximum error of 9.32\%.
Furthermore, the ROM achieves an approximate 142-fold speedup in online computations compared to the full-order model CFD simulation. 
These OpInf ROMs may be used for fast and accurate predictions of the purging flow in the PECVD chamber, which could facilitate effective particle contamination control in semiconductor manufacturing.  
\end{abstract}

\begin{keyword}
    Model order reduction \sep scientific machine learning \sep parametric Operator Inference \sep semiconductor manufacturing \sep Plasma-enhanced Chemical Vapor Deposition \sep particle contamination control
\end{keyword}

\maketitle

\section{Introduction}
\begin{figure}[!t]
\centering
\includegraphics[width=3.4in]{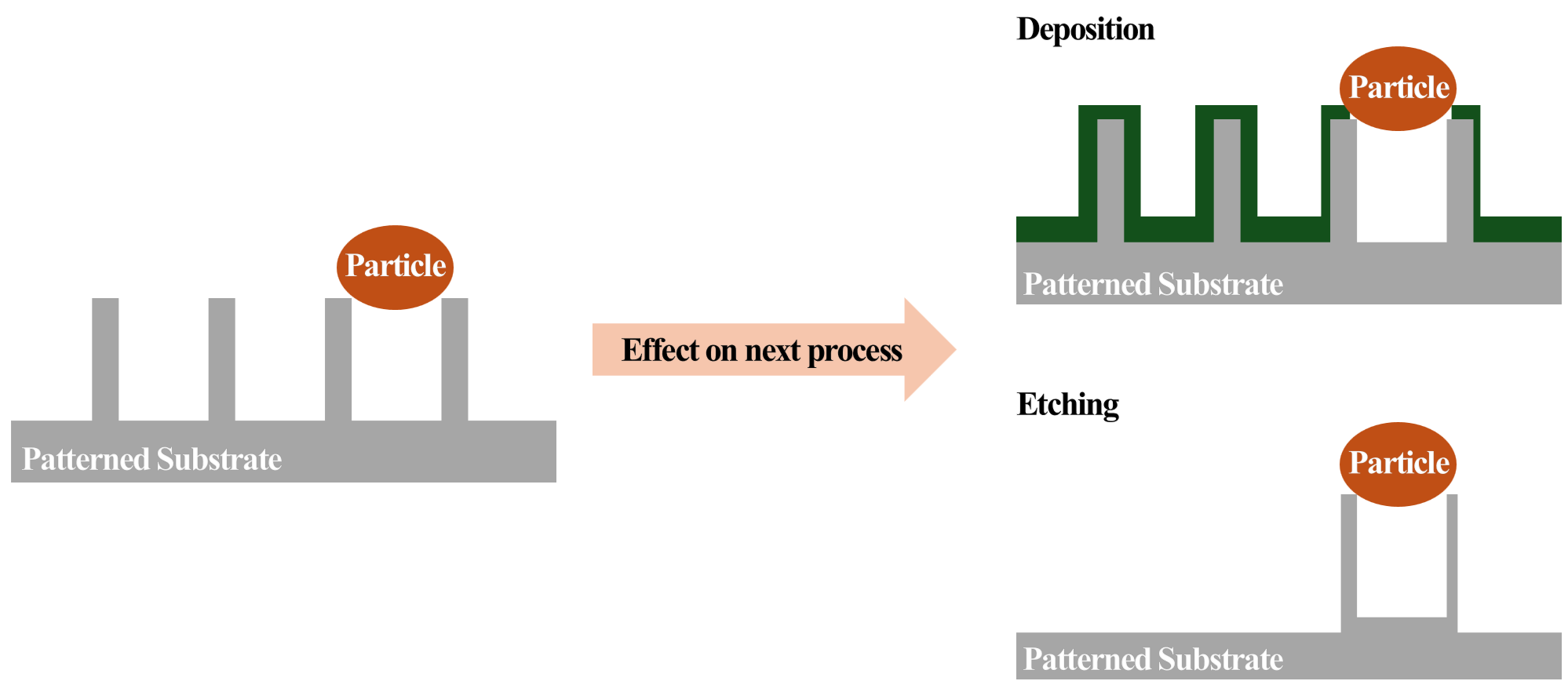}
\caption{Illustration of the impact of particle contamination on subsequent wafer manufacturing processes.}
\label{particle_contamination}
\end{figure}

Semiconductors have complex structures that are manufactured through a series of sophisticated processes, such as photolithography for pattern transfer, deposition for thin layer addition, etching for material removal, and doping to modify electrical properties. 
Across many fabrication steps, particle contamination has long been a major challenge in semiconductor manufacturing. If particles, which typically range in size from a few nanometers to several micrometers, land on the wafer surface, they disrupt the formation of the desired semiconductor structure, significantly reducing both yield and product quality. Fig.~\ref{particle_contamination} provides an illustration of the impact of particle contamination on subsequent processes in wafer manufacturing. It shows that particles settling on a patterned surface interfere with critical fabrication steps, such as deposition, where thin films are formed on the substrate, and etching, which selectively removes material to create a desired pattern. Specifically, particle-covered regions prevent material accumulation in the deposition process, while these areas remain unaffected in the etching process. 
Particle contamination in semiconductor manufacturing can originate from various sources, such as wafer transfer, equipment wear, manufacturing processes, airborne particles in the environment, etc. Among these, the particles generated during the semiconductor manufacturing process itself are the most unavoidable and high-risk sources.
This is because the wafer surface must be directly exposed to an environment where various reactions occur to form the desired structures. Therefore, it is important to ensure that particles produced during the manufacturing process do not contaminate the wafer surface and are removed from the chamber. \par
The plasma-enhanced chemical vapor deposition~(PECVD) process is a widely used semiconductor fabrication technique for depositing thin films onto a wafer surface. This process involves using plasma to enhance the chemical reaction of gaseous precursors, which allows the deposition of material at lower temperatures compared to traditional CVD. 
Especially, particles generated during the PECVD process tend to remain trapped at a certain height above the wafer due to the influence of the electric field while the plasma is active~\cite{Particle_Hiroyuki2009}. When the plasma is deactivated, the electric field dissipates, and the purging process---a subprocess of PECVD---removes these particles from the chamber using gas flow. Since this flow primarily drives the movement of particles, accurate analysis of the purging flow field during the purging process is crucial for controlling and preventing particle contamination in the PECVD process. \par

Despite the growing emphasis on particle contamination control, completely achieving this remains challenging.
From a technical standpoint, this is partially because semiconductor design rules continue to shrink, meaning that the critical dimensions of semiconductor devices, such as line widths and spacing, are decreasing. As feature sizes become smaller, particles that were previously negligible can now cause defects, increasing contamination sensitivity.
Also, the adoption of new materials and more intricate and stringent fabrication techniques has further intensified the challenges of contamination control.
From an industrial perspective, as the semiconductor industry advances toward higher efficiency and cost-effectiveness, rigorous particle contamination control has become imperative to maintain technological competitiveness.
Thus, ensuring effective particle contamination control is essential to sustain high wafer yields, preserve product quality, and enhance manufacturing productivity. This ultimately allows for the timely production and market availability of semiconductor devices. \par

To accurately identify and address particle contamination mechanisms, both experimental and numerical analyses are essential.
Over the years, many studies have leveraged experimental analysis to improve particle contamination control, focusing on contamination mechanisms~\cite{Process_Moriya2005, Process_Jeong2024, Process_Munsu2015, ku2024minimization}, particle detection during manufacturing~\cite{park2018kernel, Detecting_JianGang2021}, and particle removal~\cite{CLN_Bakhtari2006, CLN_Harald2013}.
Although experimental approaches offer direct and immediate insights into particle contamination control, their feasibility is often limited in real industrial settings due to the highly sensitive operating environment. 
Numerical simulation-based analysis, on the other hand, allows for a wider range of interrogations of the process and can provide detailed and comprehensive insights into the contamination mechanisms.
In particular, high-fidelity computational fluid dynamics (CFD) simulations have been widely used to analyze the fluid flow dynamics within manufacturing equipment, which plays a crucial role in governing particle behavior~\cite{CFD_Hojun2021, CFD_Sejin2010}.
However, the substantial computational cost of CFD simulations poses a challenge for their application in semiconductor manufacturing, where rapid problem solving is essential.
Recently, an increasing number of studies have developed machine learning models for semiconductor manufacturing by utilizing both experimental and simulated data~\cite{huang2023survey, 10601235, ding2019machine, ding2021machine, jin2024machine, tsunooka2018high}. 
These methods either uncover hidden patterns and provide insights for informed decision-making or minimize the turnaround time for analysis to enable faster decision-making.
Specifically in the context of particle control, the authors in~\cite{ML_Feng2023} develop machine learning models trained on experimental data. However, the development of machine learning models trained on numerical simulation data for particle control remains limited. \par
Our objective is to rapidly and accurately predict the flow field within the PECVD chamber during the purging process, under varying process parameters. This enables effective prediction of particle behavior, which is conducive to minimizing contamination.
In this paper, we aim to mitigate the computational bottleneck of conventional CFD simulation by proposing a data-driven reduced-order model~(ROM)---a surrogate for a full-order model (FOM)---using data obtained from CFD simulations.
In particular, we use Operator Inference (OpInf)~\cite{PEHERSTORFER2016196, kramer2024learning}---a scientific machine learning method that learns from data a low-dimensional representation of a high-dimensional system---to predict fluid dynamics within the PECVD process chamber. This nonintrusive approach builds a surrogate model without requiring access to FOM numerical operators or routines of the CFD code. 
More specifically, we build parametric ROMs to forecast the flow field within the chamber by interpolating between models with different purging process parameters.
The results demonstrate that parametric OpInf achieves high predictive accuracy for flow field prediction across varying parameters while significantly improving computational efficiency.  
    
\section{Computational model for the purging process} \label{sec:compeModel}
\subsection{Computational Domain of the PECVD Chamber} \label{ss:computeDomain} 
\begin{figure*}[!t]
\centering
\subfloat[]{\includegraphics[height=2.5in, keepaspectratio]{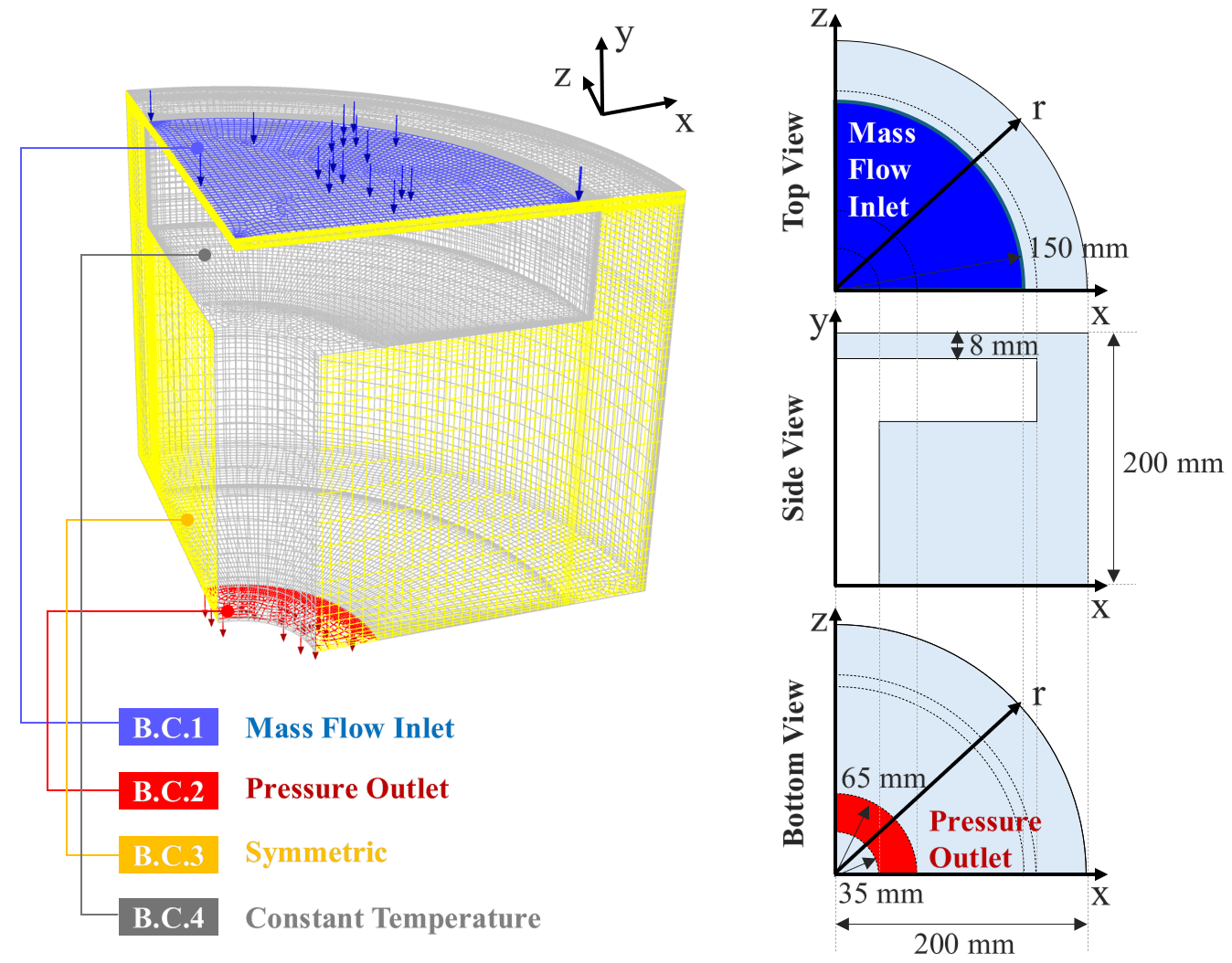}%
\label{fig:computational_domain}}
\hfil
\subfloat[]{\includegraphics[height=2.5in, keepaspectratio]{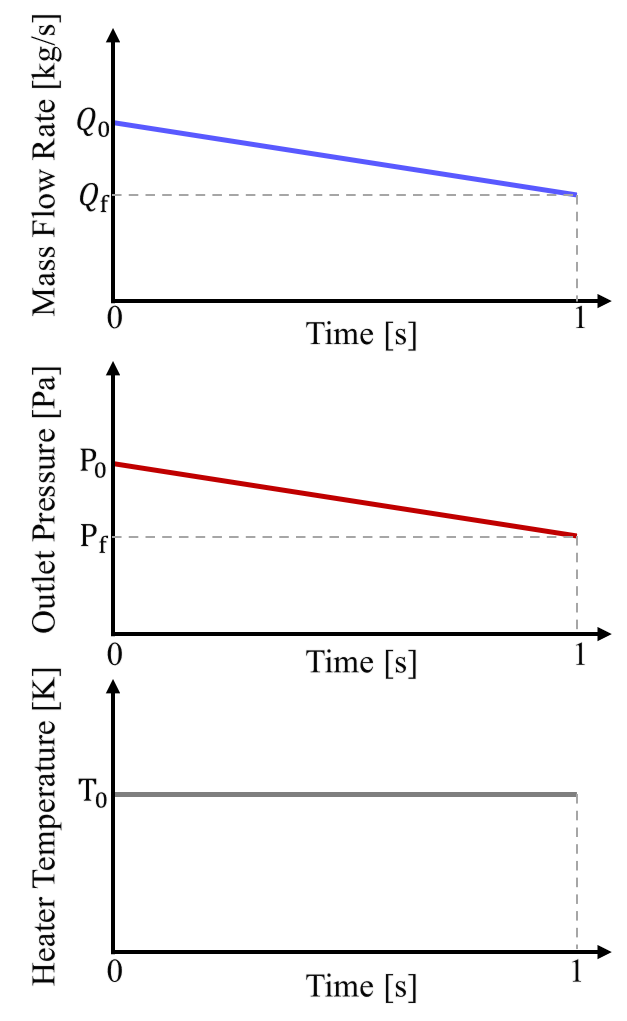}%
\label{fig:BC_profile}}
\caption{(a) Computational domain of a simplified PECVD chamber. B.C.1: mass flow inlet boundary condition, B.C.2: pressure outlet boundary condition, B.C.3: symmetric boundary condition. The side view represents the $xy$-vertical plane where $z~=~0$~mm. (b) Profiles of the mass flow inlet, pressure outlet, and heater surface temperature over time $t \in [0, 1]s$.}
\label{fig:domain_and_BC}
\end{figure*}

We simulate the purging process in a simplified chamber structure, as illustrated in Fig.~\ref{fig:computational_domain}. Given that this chamber has a cylindrical shape with a total radius of 200~mm, its physical properties and governing equations remain unchanged across symmetric planes. This allows us to model only one quarter of the chamber, reducing computational cost in the FOM while still capturing the essential flow behavior. We, thus, use symmetric boundary conditions.
The top surface, with a radius of 150~mm, serves as the argon (Ar) gas mass flow inlet, while the bottom surface, which covers the radius range of 35 to 65~mm, acts as the pressure outlet. A wafer heater, positioned 8 mm below the top surface, has a radius of 160~mm and supports the wafer during the purging process. 
We generate a total of 249,262 mesh elements that consist of 232,304 8-node hexahedral elements and 1,448 6-node wedge elements.

\begin{figure}[!t]
\centering
\includegraphics[height=1.5in]{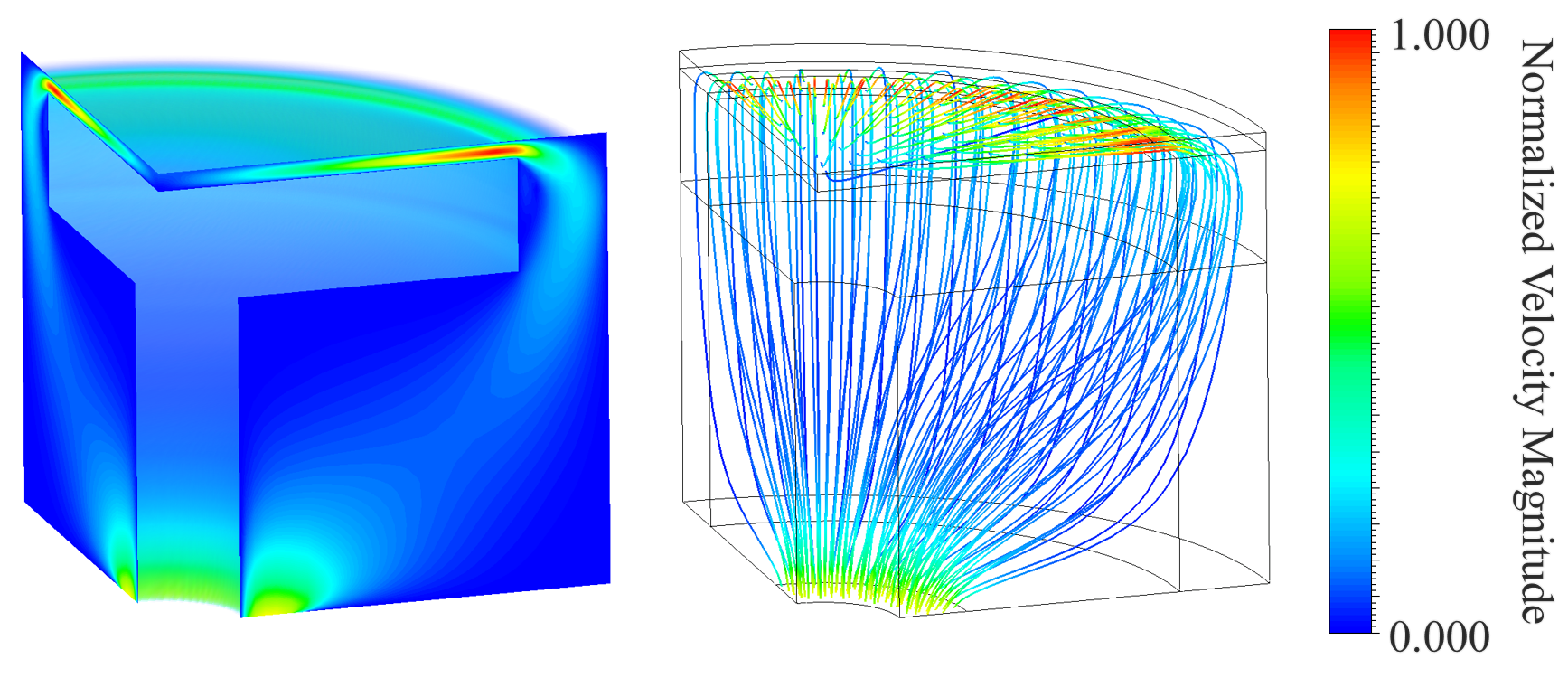}
\caption{Left: velocity field, Right: streamlines at the onset of the purging process, initiated upon plasma deactivation.}
\label{fig:FOM_results_initial}
\end{figure}

\subsection{CFD Simulation of the Purging Process} \label{ss:compMethod}
The purging process within the PECVD chamber occurs after plasma deactivation and the completion of deposition. During purging, chemical reactions are negligible, and particle transport is assumed to be driven mainly by hydrodynamic drag from the bulk gas flow.
The low operating pressures and flow rates typical of purging conditions yield Reynolds numbers well within the laminar regime. The resulting flow field, characterized by pressure, temperature, and velocity distributions, therefore provides the critical insights necessary to inform particle transport in practice.
We simulate the purging process within the PECVD chamber using the commercial software ANSYS Fluent 2023 R1, which solves the governing equations for mass, momentum, and energy using the finite volume method. The flow is assumed to be laminar, and the energy equation is included to capture the effects of heat transfer.
When the purging process begins, the Ar mass flow inlet and the pressure outlet boundary conditions change linearly over time, as shown in Fig.~\ref{fig:BC_profile}. These time-dependent variations are implemented using the boundary condition profile in ANSYS Fluent, while the heater surface is maintained at a constant temperature to ensure stable thermal conditions.
To predict the internal flow field within the chamber depending on variations of the inlet mass flow and the outlet pressure, we conduct FOM simulations across a range of two parameter values, denoted as $\mu_q$ and $\mu_p$. These parameters define the scaling of the final inlet mass flow rate and outlet pressure relative to their initial values. The final inlet Ar mass flow rate (see Fig.~\ref{fig:BC_profile}) is given by $Q_f = \mu_q Q_0$, and the final outlet pressure by $P_f = \mu_p P_0$, where $Q_0$ and $P_0$ are the initial inlet Ar mass flow rate and outlet pressure, respectively. Fig.~\ref{fig:FOM_results_initial} shows the initial flow field in the PECVD chamber. When the purging process begins, steady airflow is formed, with the Ar gas entering through the upper inlet and exiting through the bottom exhaust.

\subsection{Parameter-dependent Full-order Model}    \label{ss:problem_setup}
Consider a parameter vector $\boldsymbol{\mu} \in \cM \subset \real^{d_p}$ of dimension $d_p \in \nat$. Given the two parameters $\mu_q$ and $\mu_p$ discussed in~\Cref{ss:compMethod}, our parametric dimension is $d_p=2$.
The semi-discretization of the nonlinear governing partial differential equations (PDEs) results in an $N$-dimensional system of nonlinear ordinary differential equations (ODEs) written generically as 
\begin{align}
    \frac{\mathrm{d} \bs(t;\boldsymbol{\mu})}{\mathrm{d} t} = \bf(\bs(t;\boldsymbol{\mu}), \bu(t;\boldsymbol{\mu})), \hspace{2em} \bs(t_0;\boldsymbol{\mu})=\bs_0(\boldsymbol{\mu}),
\label{eq:FOM}
\end{align}
which describes the time evolution of the discretized full state vector $\bs(t;\boldsymbol{\mu}) = [s_1(t;\boldsymbol{\mu}),\ldots, s_{N}(t;\boldsymbol{\mu})]^\top \in \real^{N}$, where $N$ is large.
Here, the $m$ inputs $\bu(t;\boldsymbol{\mu})$ result from the time-dependent boundary conditions, and the nonlinear function $\bf:\real^N \times [t_0, t_f] \times \cM \times \real^{m} \rightarrow \real^N$ maps the parametrized full state $\bs(t;\boldsymbol{\mu})$ and the input $\bu(t;\boldsymbol{\mu})$ to the time derivative of the full state.

\subsection{Full-order Model Data from CFD Simulation} \label{ss:FOM_details}
The high-fidelity data consists of $n_x = 249,262$ semi-discretized spatial elements and $n_v = 5$ variables (pressure, temperature, three velocity components) at each time step and each combination of parameters. Thus, $N = n_x \cdot n_v = 1{,}246{,}310$ becomes the degrees of freedom of the high-fidelity model. The simulation is conducted over $T = 1s$, with a constant time step $\delta = 0.005s$, resulting in $K = 200$ time steps during the purging process.
We generate a total of $25$ FOM solutions for five distinct values for both the flow rate parameter~$\mu_q$ and the pressure parameter~$\mu_p$. From the FOM data, $d = 9$ solutions are used for training, while 16 are used for testing. \par

Given a set of training parameter vectors $\{{\boldsymbol{\mu}_\ell}\}_{\ell=1}^{d}$, where $ \boldsymbol{\mu}_\ell \in \cM$ for each $\ell=1,\ldots,d$, we collect the state trajectories for all $d$ parameter instances, and store them in the global data matrix
\begin{align*}
    \bS = \left[ \bS(\boldsymbol{\mu}_1) \hspace{0.7em} \bS(\boldsymbol{\mu}_2) \hspace{0.7em} \cdots \hspace{0.7em} \bS(\boldsymbol{\mu}_d) \right] \in \real^{N \times dK},
\end{align*}
where each $\bS(\boldsymbol{\mu}_\ell) = [\bs(t_0;\boldsymbol{\mu}_\ell), \bs(t_1;\boldsymbol{\mu}_\ell), \ldots, \bs(t_{K-1}; \boldsymbol{\mu}_\ell)] \in \real^{N \times K}$. For ease of notation, we assume that each simulation requires the same number of time steps, but the framework directly carries over to the case where this does not hold.

\section{Parametric Operator Inference for nonintrusive learning of the PECVD purging~process} \label{sec:OpInf}
Since we simulate the purging process via ANSYS Fluent, extracting the high-fidelity operators and functions of the semi-discretized PDE is not feasible.
For such a setting, Operator Inference~(OpInf) has emerged as an efficient scientific machine learning method that constructs predictive polynomial ROMs of high-dimensional dynamical systems~\cite{PEHERSTORFER2016196, kramer2024learning}. 
We apply parametric OpInf to interpolate across multiple datasets generated from simulations with varying parameter values, where the parameters $\mu_q$ and $\mu_p$ represent the scaling factors between the initial and final values of the time-varying boundary conditions.
Previous studies have incorporated parametric dependencies into the OpInf framework using different approaches:~\cite{mcquarrie2023nonintrusive} uses an affine formulation to embed the parametric structure of the governing equations directly into the regression problem,~\cite{yildiz2021learning} assumes an affine dependence of the reduced operators on the parameters of interest, and~\cite{farcas2023parametric} utilizes a Taylor series expansion to represent the parametric dependence of the reduced operators.
However, since the problem considered herein does not exhibit affine dependence, we adopt an interpolatory regression approach, as we discuss~next.

\subsection{Learning Parametrized ROMs via Interpolatory Regression} \label{ss:parametric_interpolation}
\subsubsection{Dimensionality reduction} \label{sss:dimension_reduction}
To reduce the dimensionality of the data, we use proper orthogonal decomposition (POD), which determines its basis based on a reduced singular value decomposition (SVD) of the global data matrix $\bS$ as
\begin{align}
    \bS = \bV \mathbf{\Sigma} \bW^\top.
\label{eq:SVD}
\end{align}
Here, the columns of $\bV \in \mathbb{R}^{N \times d K}$ and  $\bW \in \real^{d K \times d K}$ are the left and right singular vectors of $\bS$, respectively, and $\mathbf{\Sigma} \in \mathbb{R}^{d K \times d  K}$ is a diagonal matrix with the singular values of $\bS$, denoted as $\sigma_1 \geq \sigma_2 \geq \cdots \geq \sigma_{d K} \geq 0$.
The global POD basis $\bV_r \in \real^{N \times r}$ consists of the leading $r$ left singular vectors in $\bV$, which correspond to the $r$ largest singular values of $\bS$, where $r \ll N$. This is an orthogonal basis, such that $\bV_r^\top \bV_r = \bI_r$ where $\bI_r \in \real^{r \times r}$ is the identity matrix.
We project the high-dimensional data $\bS$ onto the linear subspace spanned by the columns of $\bV_r$, yielding the low-dimensional data matrix
\begin{align*}
    \widehat{\bS} = \bV_r^\top \bS = \left[ \widehat{\bs}(t_0; \boldsymbol{\mu}) \hspace{0.5em} \widehat{\bs}(t_1; \boldsymbol{\mu}) \hspace{0.5em} \cdots \hspace{0.5em} \widehat{\bs}(t_{dK-1}; \boldsymbol{\mu}) \right] \in \real^{r \times dK}.
\end{align*}
Here, each column $\widehat{\bs}(t_k;\boldsymbol{\mu}) \in \real^{r}$, $k=0,1,\ldots,dK-1$, serves as a $\textit{reduced state}$ in the low-dimensional representation of the high-dimensional FOM~\eqref{eq:FOM}.

\subsubsection{Parametrized Quadratic ROM} \label{sss:quadratic_ROM}
The OpInf framework learns a ROM in polynomial form, and choosing a quadratic polynomial is often a good trade-off between model expressiveness and computational runtime for the ROM simulation \cite{McQuarrie_2021, Swischuk_2020, qian2022reduced, farcas2023parametric, Zastrow2023, 0x003d183f, Rocha2023}. 
Thus, we opt to learn a parameter-dependent quadratic ROM
\begin{align*}
    \frac{\mathrm{d}\widehat{\bs}(t;\boldsymbol{\mu})}{\mathrm{d}t} &= \widehat{\bc}(\boldsymbol{\mu}) + \widehat{\bA}(\boldsymbol{\mu}) \widehat{\bs}(t;\boldsymbol{\mu}) + \widehat{\bH}(\boldsymbol{\mu}) \left( \widehat{\bs}(t;\boldsymbol{\mu}) \otimes \widehat{\bs}(t;\boldsymbol{\mu}) \right) \\ &\quad + \widehat{\bB}(\boldsymbol{\mu}) \bu(t;\boldsymbol{\mu}).
\end{align*}
Here, the operator $\otimes$ indicates a compact Kronecker product that calculates the elementwise multiplication of two vectors while removing redundant terms. For example, when $\ba = [a_1, a_2, a_3]$, the compact Kronecker product is computed as $\ba \otimes \ba = [a_1^2, a_1a_2, a_1a_3, a_2^2, a_2a_3, a_3^2]$. 
The $\widehat{\bc}(\boldsymbol{\mu}) \in \real^r$ is a reduced constant vector, $\widehat{\bA}(\boldsymbol{\mu}) \in \real^{r \times r}$ is a reduced linear operator, $\widehat{\bH}(\boldsymbol{\mu}) \in \mathbb{R}^{r \times r(r+1)/2}$ is a reduced quadratic operator, and $\widehat{\bB}(\boldsymbol{\mu}) \in \real^{r \times m}$ is a reduced input operator.

\subsubsection{ROM Learning via Linear Least-squares Regression}   \label{sss:least_squares_regression}
To learn the reduced operators which are parametrized by training parameters $\boldsymbol{\mu}_\ell \in \cM$, for $\ell = 1,\ldots, d$, OpInf solves the linear least-squares problem 
\begin{equation}
\begin{aligned}
    &\min_{\widehat{\bc}_\ell, \widehat{\bA}_\ell, \widehat{\bH}_\ell, \widehat{\bB}_\ell} 
    \sum_{\ell=1}^{d}
    \sum_{k=0}^{dK-1} 
    \Big\|
        \widehat{\bc}_\ell + \widehat{\bA}_\ell \widehat{\bs}(t_k; \boldsymbol{\mu}_\ell) + \\
        & \widehat{\bH}_\ell 
        \big( 
            \widehat{\bs}(t_k; \boldsymbol{\mu}_\ell) 
            \otimes 
            \widehat{\bs}(t_k; \boldsymbol{\mu}_\ell) 
        \big) 
        + \widehat{\bB}_\ell \bu(t_k; \boldsymbol{\mu}_\ell) 
        - \dot{\widehat{\bs}}(t_k; \boldsymbol{\mu}_\ell) 
    \Big\|_2^2,
\end{aligned}
\label{eq:parametric_least_squares_problem}
\end{equation}
where $\widehat{\bc}_\ell = \widehat{\bc}(\boldsymbol{\mu}_\ell)$, $\widehat{\bA}_\ell=\widehat{\bA}(\boldsymbol{\mu}_\ell)$, $\widehat{\bH}_\ell=\widehat{\bH}(\boldsymbol{\mu}_\ell)$, and $\widehat{\bB}_\ell=\widehat{\bB}(\boldsymbol{\mu}_\ell)$.
In the previous equation, $\dot{\widehat{\bs}}(t_k;\boldsymbol{\mu}_\ell)$ represents the time derivative of the reduced state $\widehat{\bs}(t;\boldsymbol{\mu}_\ell)$ at time $t = t_k$, and it can be computed using a suitable time derivative approximation scheme.~\Cref{eq:parametric_least_squares_problem} can be written compactly as
\begin{align}
    \min_{\widehat{\bO}_1, \ldots, \widehat{\bO}_d} & \underbrace{ \sum_{\ell=1}^{d} \left\| \bD_\ell \widehat{\bO}_\ell^\top -  \dot{\widehat{\bS}}_\ell^\top \right\|_F^2 }_{L(\widehat{\bO}_\ell)},
\label{eq:parametrized_regression_matrix_form}
\end{align}
where the reduced operators matrix~$\widehat{\bO}_\ell=\widehat{\bO}(\boldsymbol{\mu}_\ell)=[ \widehat{\bc}_\ell \hspace{0.5em} \widehat{\bA}_\ell \hspace{0.5em} \widehat{\bH}_\ell \hspace{0.5em} \widehat{\bB}_\ell ] \in \real^{r \times d(r,m)}$ and $d(r,m) = 1~+~r + r(r+1)/2 + m$.
The data matrix $\bD_\ell = [\mathbf{1}_{dK} \hspace{0.5em} \widehat{\bS}_\ell^\top \hspace{0.5em} (\widehat{\bS}_\ell \otimes \widehat{\bS}_\ell)^\top \hspace{0.5em} \bU_\ell^\top] \in \real^{dK \times d(r,m)}$, where $\mathbf{1}_{dK} \in \mathbb{R}^{dK}$ is a column vector of length $dK$ with all entries set to one, $\widehat{\bS}_\ell = \widehat{\bS}(\boldsymbol{\mu}_\ell) = [\widehat{\bs}(t_0; \boldsymbol{\mu}_\ell) \cdots \widehat{\bs}(t_{dK-1}; \boldsymbol{\mu}_\ell)] \in \real^{r \times dK}$, and $\bU_\ell = \bU(\boldsymbol{\mu}_\ell) = [\bu(t_0;\boldsymbol{\mu}_\ell) \cdots \bu(t_{dK-1};\boldsymbol{\mu}_\ell)] \in \real^{m \times dK}$.
Also, $\dot{\widehat{\bS}}_\ell \in \real^{r \times dK}$ is a matrix whose columns are $\dot{\widehat{\bs}}(t_k;\boldsymbol{\mu}_\ell)$.
We then rewrite the cost function $L(\widehat{\bO}_\ell)$ in~\eqref{eq:parametrized_regression_matrix_form} as
\begin{equation*}
    L(\widehat{\bO}_\ell) = \sum_{\ell=1}^{d} \underbrace{\text{tr} \left( \bD_\ell \widehat{\bO}_\ell^\top -  \dot{\widehat{\bS}}_\ell^\top \right)^\top \left( \bD_\ell \widehat{\bO}_\ell^\top -  \dot{\widehat{\bS}}_\ell^\top \right) }_{f(\widehat{\bO}_\ell)}.
\end{equation*}    
To minimize $L(\widehat{\bO}_\ell)$, we apply the first-order optimality condition, which involves taking the gradient of $L(\widehat{\bO}_\ell)$ with respect to each $\widehat{\bO}_\ell$ and setting it to zero.
The gradient of each term in the cost function with respect to $\widehat{\bO}_\ell$ is simplified as
\begin{equation*}
\begin{aligned}
    \nabla_{\widehat{\bO}_\ell} f(\widehat{\bO}_\ell) 
    &= \nabla_{\widehat{\bO}_\ell} 
    \text{tr} \big( 
        \widehat{\bO}_\ell \bD_\ell^\top \bD_\ell \widehat{\bO}_\ell^\top 
        - \widehat{\bO}_\ell \bD_\ell^\top \dot{\widehat{\bS}}_\ell^\top 
        - \dot{\widehat{\bS}}_\ell \bD_\ell \widehat{\bO}_\ell^\top \\
        &\quad+ \dot{\widehat{\bS}}_\ell \dot{\widehat{\bS}}_\ell^\top  
    \big) \\
    &= \nabla_{\widehat{\bO}_\ell} 
    \Big[ 
        \text{tr} \big( \widehat{\bO}_\ell \bD_\ell^\top \bD_\ell \widehat{\bO}_\ell^\top \big) 
        - \text{tr} \big( \widehat{\bO}_\ell \bD_\ell^\top \dot{\widehat{\bS}}_\ell^\top \big) \\
        &\quad - \text{tr} \big( \dot{\widehat{\bS}}_\ell \bD_\ell \widehat{\bO}_\ell^\top \big) 
        + \text{tr} \big( \dot{\widehat{\bS}}_\ell \dot{\widehat{\bS}}_\ell^\top \big) 
    \Big] \\
    &= \widehat{\bO}_\ell 
    \big( \bD_\ell^\top \bD_\ell + \bD_\ell^\top \bD_\ell \big) 
    - \big( \bD_\ell^\top \dot{\widehat{\bS}}_\ell^\top \big)^\top 
    - \dot{\widehat{\bS}}_\ell \bD_\ell \\ 
    &= 2 \big( \widehat{\bO}_\ell \bD_\ell^\top \bD_\ell - \dot{\widehat{\bS}}_\ell \bD_\ell \big).
\end{aligned}
\end{equation*}
Thus, the first-order optimality condition $\nabla_{\widehat{\bO}_\ell} L(\widehat{\bO}_\ell) = 0$ is solved by 
\begin{equation}
    \bD_\ell^\top \bD_\ell \widehat{\bO}_\ell^\top = \bD_\ell^\top \dot{\widehat{\bS}}_\ell^\top ,
\label{eq:normal_equation}
\end{equation}
for each $\widehat{\bO}_\ell$, where $\ell = 1,\ldots,d$.
In other words,~\eqref{eq:parametrized_regression_matrix_form} decouples into $d$ individual normal equations~\eqref{eq:normal_equation}.

\subsubsection{ROM Interpolation}   \label{sss:ROM_interpolation}
Once the reduced operators $\widehat{\bO}_1, \ldots, \widehat{\bO}_d$ are learned for all $d$ parameters~$\boldsymbol{\mu}_1, \ldots, \boldsymbol{\mu}_d$, the reduced model for $\boldsymbol{\mu} \in \cM$ is derived via elementwise interpolation~\cite{peherstorfer2020sampling} as 
\begin{align*}
    \widehat{\bO}(\boldsymbol{\mu}) = \textnormal{interpolate}((\boldsymbol{\mu}_1, \widehat{\bO}_1), \ldots, (\boldsymbol{\mu}_d, \widehat{\bO}_d) ; \boldsymbol{\mu}) .
\end{align*}
For the implementation of elementwise interpolation between the parameter vectors $\boldsymbol{\mu}_\ell$ and the reduced operators $\widehat{\bO}_\ell$, we use the~\texttt{LinearNDInterpolator} module~\cite{10.1145/235815.235821} in the $\mathtt{scipy}$ package~\cite{virtanen2020scipy} in Python.
This interpolator constructs the ROM with reduced operators $\widehat{\bc}^*$, $\widehat{\bA}^*$, $\widehat{\bH}^*$, and $\widehat{\bB}^*$ for a given parameter vector $\boldsymbol{\mu}^* = (\mu_q^*, \mu_p^*) \in \cM$, enabling state prediction for the new parameter unseen during the training. For the implementation of the parametric ROM learning and state predictions, we use the~\texttt{opinf} Python package version~0.5.11~\cite{OpInfGitHub}.

\subsection{Data Transformation for Multiscale Data}  \label{ss:data_transformation}
The purging flow behavior in the PECVD chamber is multiphysics and multiscale, involving variables with different orders of magnitude. This poses a substantial challenge in model training, as variables with larger magnitudes can dominate those with smaller scales. In our case, the original variables, temperature and velocities, differ by up to four orders of magnitude.
To mitigate numerical issues, we pre-process the training data by scaling all variables to comparable orders of magnitude. Specifically, we apply mean subtraction to the pressure and temperature fields, and scale all variables to the range of~[$-1,1$]. This results in the scaled data matrix~$\bX \in \mathbb{R}^{1{,}246{,}310 \times 1{,}800}$, which is used to construct the OpInf ROM. 

\subsection{OpInf ROM Learning via Regularization}  \label{ss:regularization}
In the model learning process~\eqref{eq:parametrized_regression_matrix_form}, we use Tikhonov regularization to avoid overfitting, which arises from potential noise, under-resolution of the data, the truncated global POD modes that result in unresolved system dynamics, and mis-specification of the ROMs as fully quadratic~\cite{McQuarrie_2021}.
Thus, the individual regression problem for each parameter $\boldsymbol{\mu}_\ell$ in~\eqref{eq:parametric_least_squares_problem} becomes
\begin{align}
    \min_{\widehat{\bO}_\ell} & \left\| \bD_\ell \widehat{\bO}_\ell^\top -  \dot{\widehat{\bX}}_\ell^\top \right\|_F^2 + \left \| \mathbf{\Lambda} \widehat{\bO}_\ell^\top \right \|_F^2,
\label{eq:optimization_regularized}
\end{align}
where $\mathbf{\Lambda} = \textnormal{diag}(\lambda_1, \lambda_1 \bI_{r}, \lambda_2 \bI_{r(r+1)/2}, \lambda_3 \bI_{m}) \in \mathbb{R}^{d(r,m) \times d(r,m)}$ is the diagonal matrix used for regularization. Here, $\bI_{r}$, $\bI_{r(r+1)/2}$, and $\bI_m$ are identity matrices with dimensions $r$, $r(r+1)/2$, and $m$, respectively.
Note that the Kronecker product $\otimes$ introduces scaling differences between the entries of the quadratic operator $\widehat{\bH}$, and those in the constant vector $\widehat{\bc}$ and the linear operator $\widehat{\bA}$, when subjected to a single regularization parameter. It is therefore best practice to use different regularization parameters $\lambda_1 > 0$, $\lambda_2 > 0$, and $\lambda_3 > 0$, such that $\lambda_1$ penalizes $\widehat{\bc}$ and $\widehat{\bA}$, $\lambda_2$ penalizes $\widehat{\bH}$, and $\lambda_3$ penalizes $\widehat{\bB}$~\cite{McQuarrie_2021}.
To determine the hyperparameters $\lambda_1$, $\lambda_2$, and $\lambda_3$ for regularizing the least-squares problem in~\eqref{eq:optimization_regularized}, we generate a uniformly spaced parameter grid on a logarithmic scale, where ($\lambda_1, \lambda_2, \lambda_3) \in [10^{-3}, 10^{3}] \times [10^{0} \times 10^{6}] \times [10^{-3} \times 10^3]$, with each dimension discretized into seven values. 
We find the regularization hyperparameters that minimize the relative state error on the training data. However, since the small magnitudes of the scaled data $\bX$ can skew small errors, we compute the relative state error with the unscaled variables. 
To prevent variables with higher orders of magnitude from dominating those with smaller magnitudes in their original scales, we compute the relative error for each variable separately as
\begin{equation}
    \cE^{\text{rel},i} = \frac{\left\|\bS^{\text{FOM},i} - \bS^{\text{ROM},i}\right\|_F}{\left\|\bS^{\text{FOM},i}\right\|_F},
\label{eq:state_error}
\end{equation}
where $\bS^{.,i} \in \real^{n_x \times K}$ is the $i$th variable ($i=1$ for pressure, $i=2$ for temperature, $i=3$ for y-velocity, and $i=4$ for radial velocity). Note that since $v_x$ and $v_z$ are axisymmetric, we combine their effects on the total relative state error by computing a single radial velocity: $v_r=\sqrt{v_x^2+v_z^2}$. For ease of error comparison, we compute the average of the four relative errors, as $\frac{1}{4} \sum_{i=1}^{4} \cE^{\text{rel}, i}$.
The selected hyperparameters are ($\lambda_1, \lambda_2, \lambda_3$) $=$ ($10^{-3}, 10^{2}, 10^{-3}$), which minimizes the training error across all $d$ ROMs during the model learning process.

\section{Numerical Results} \label{sec:numerical_results}
\subsection{Selection of ROM Dimension}    \label{ss:basis_selection}

\begin{figure}[!t]
\centering
\hspace{-3em}\subfloat[]{\includegraphics[height=1.8in, keepaspectratio]{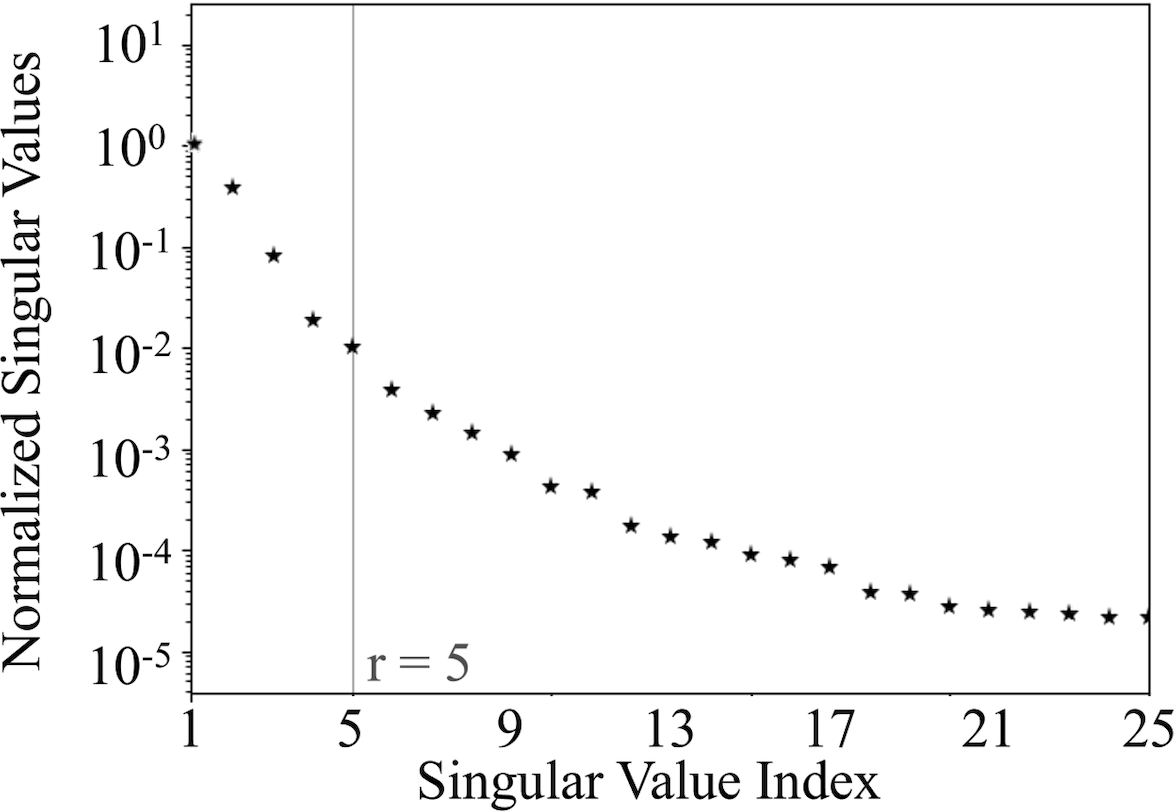}%
\label{fig:singular_value}}
\vspace{0.5em}
\subfloat[]{\includegraphics[height=1.80in, keepaspectratio]{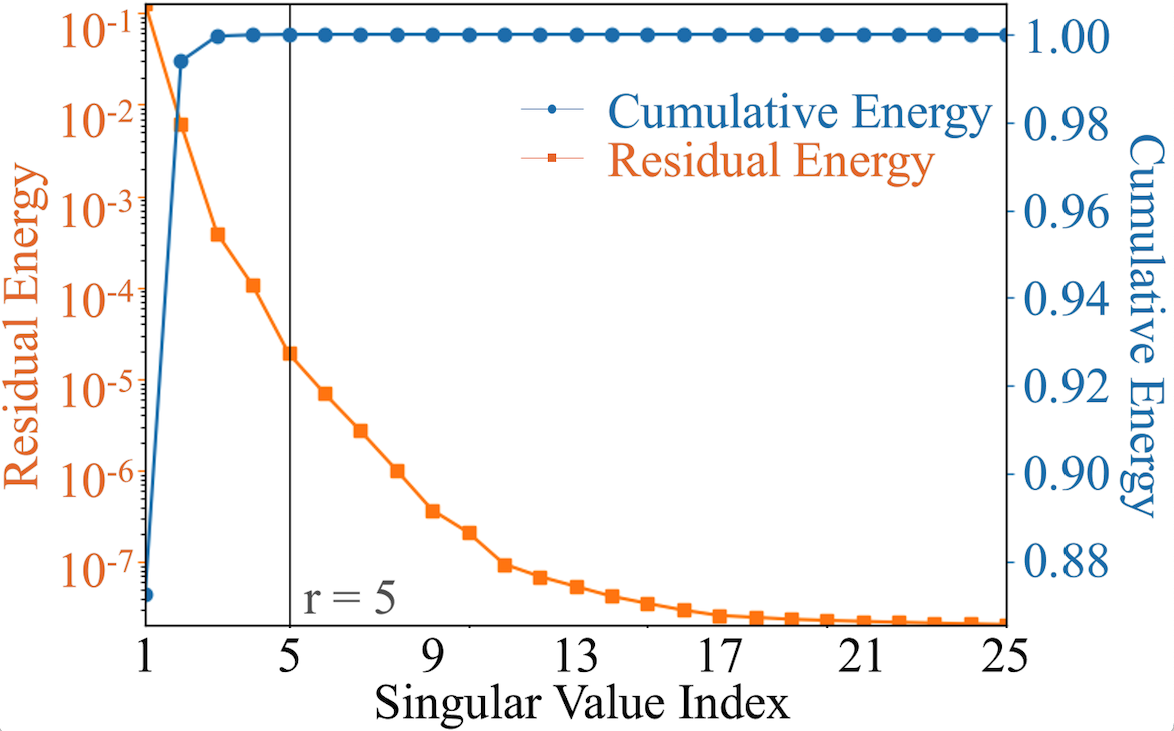}%
\label{fig:cumulative_energy}}
\caption {(a) Normalized singular value decay, (b) the cumulative energy $\cE_{\text{cum}}(r)$ and the residual energy that represents the projection error $\cE_{\text{proj}}(r)$, with respect to the reduction dimension $r$. The vertical black line illustrates the chosen ROM dimension $r=5$.}
\label{fig:singluar_value_and_cumulative_energy}
\end{figure}

We determine the ROM dimension $r$ using the energy metric, which is defined by the singular values of the data matrix. To achieve this, we first compute the SVD of the scaled training data matrix $\bX \in \real^{N \times dK}$, as described in~\eqref{eq:SVD}.
However, the deterministic SVD of $\bX$ has a complexity of \(O(N^2(dK) + (dK)^3)\), making it impractical due to the computational load that grows quadratically in~$N$. Thus, we compute a randomized SVD via the python package \texttt{Scikit-learn}~\cite{halko2011finding, pedregosa2018scikitlearnmachinelearningpython}.
From the randomized SVD of $\bX$, the cumulative energy $\cE_{\text{cum}}(r)$ is defined as the ratio of the sum of squared singular values up to rank $r$ to the total sum of squared singular values: $\cE_{\text{cum}}(r) = \sum_{\eta=1}^{r}\sigma_{\eta}^2 / \sum_{\eta=1}^{r_t}\sigma_{\eta}^{2}$. Here, $\sigma_{\eta}$ ($\eta=1,2,\dots,r_t$) are the singular values of the scaled data matrix $\bX$, and $r_t$ is a target rank when computing randomized SVD. 
This relates the low-rank dimension $r$ to the projection error of the scaled data matrix $\bX$ as follows:
\begin{equation}
    \cE_{\text{proj}}(r) = \frac{\left\| \bX - \bV_r \bV_r^{\top} \bX \right\|_F^2}{\left\| \bX \right\|_F^2} = 1 - \cE_{\text{cum}}(r),
\label{eq:energy_criterion}
\end{equation}
which quantifies the accuracy with which the linear basis $\bV_r$ reconstructs the original data $\bX$ through the linear injection of the projected data $\widehat{\bX} = \bV_r^\top \bX \in \real^{r \times dK}$.
Fig.~\ref{fig:singluar_value_and_cumulative_energy} shows the normalized singular value decay and the cumulative energy $\cE_{\text{cum}}(r)$ related to the projection energy $\cE_{\text{proj}}(r)$ in terms of the ROM dimension $r$.
The ROM dimension is set to $r=5$, which retains a cumulative energy of $99.9981\%$ (residual energy of $1.9056\times10^{-5}$), providing an accurate low-dimensional representation of the FOM dynamics while substantially reducing computational cost. 
We note that while increasing $r$ beyond 5 further reduces the residual energy, it does not necessarily improve predictive accuracy in this type of nonlinear flow application. We tested ROMs with $r=10$ and $r=15$ for a single parameter case; neither improved the relative state error~\eqref{eq:state_error} compared to $r=5$.

\subsection{Prediction Results of Parametric OpInf ROM Interpolation} \label{ss:Comparison}
\begin{figure*}[!t]
\centering
\subfloat[]{\includegraphics[height=1.6in, keepaspectratio]{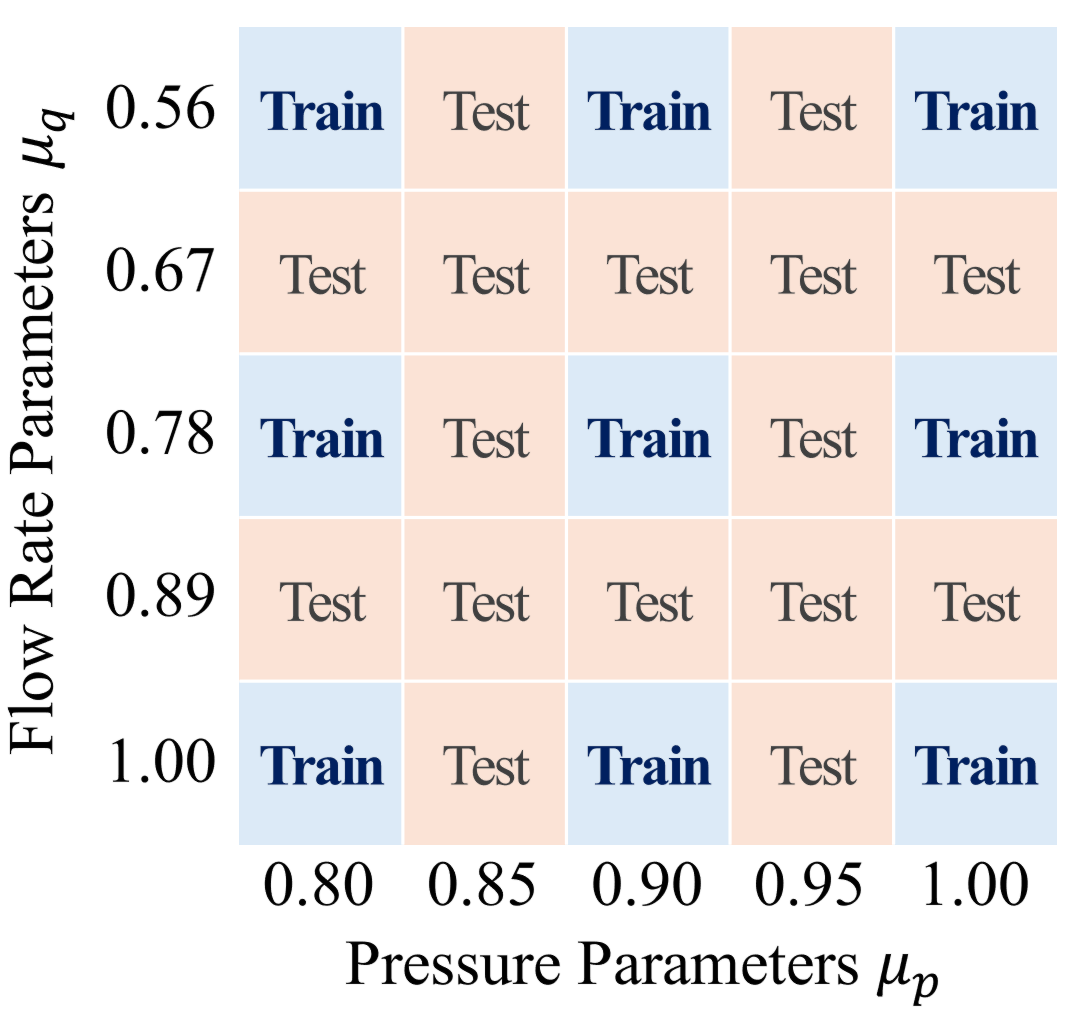}%
\label{fig:dataset_selection}}
\hfil
\subfloat[]{\includegraphics[height=1.6in, keepaspectratio]{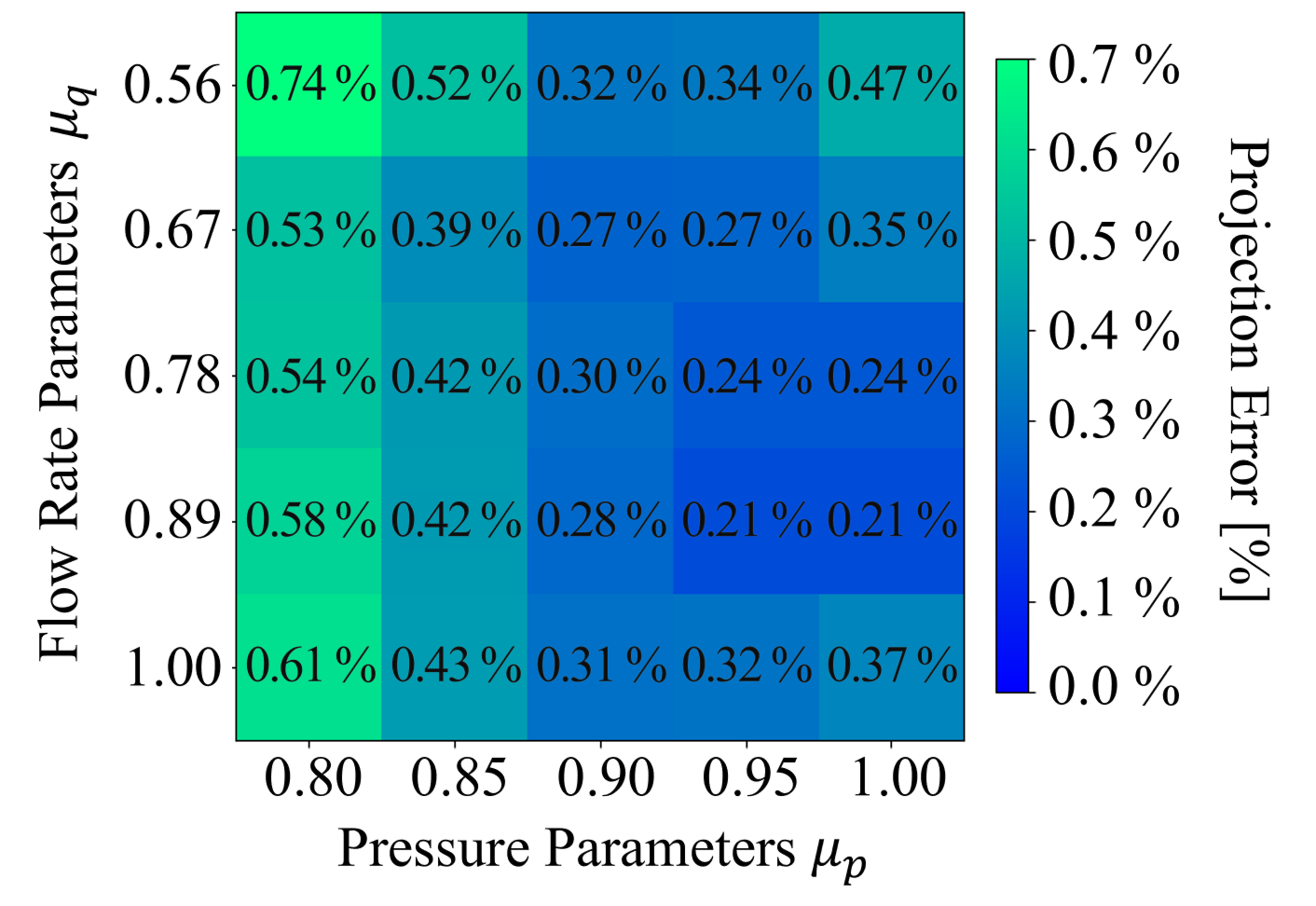}%
\label{fig:projection_error_all}}
\hfil
\subfloat[]{\includegraphics[height=1.6in,keepaspectratio]{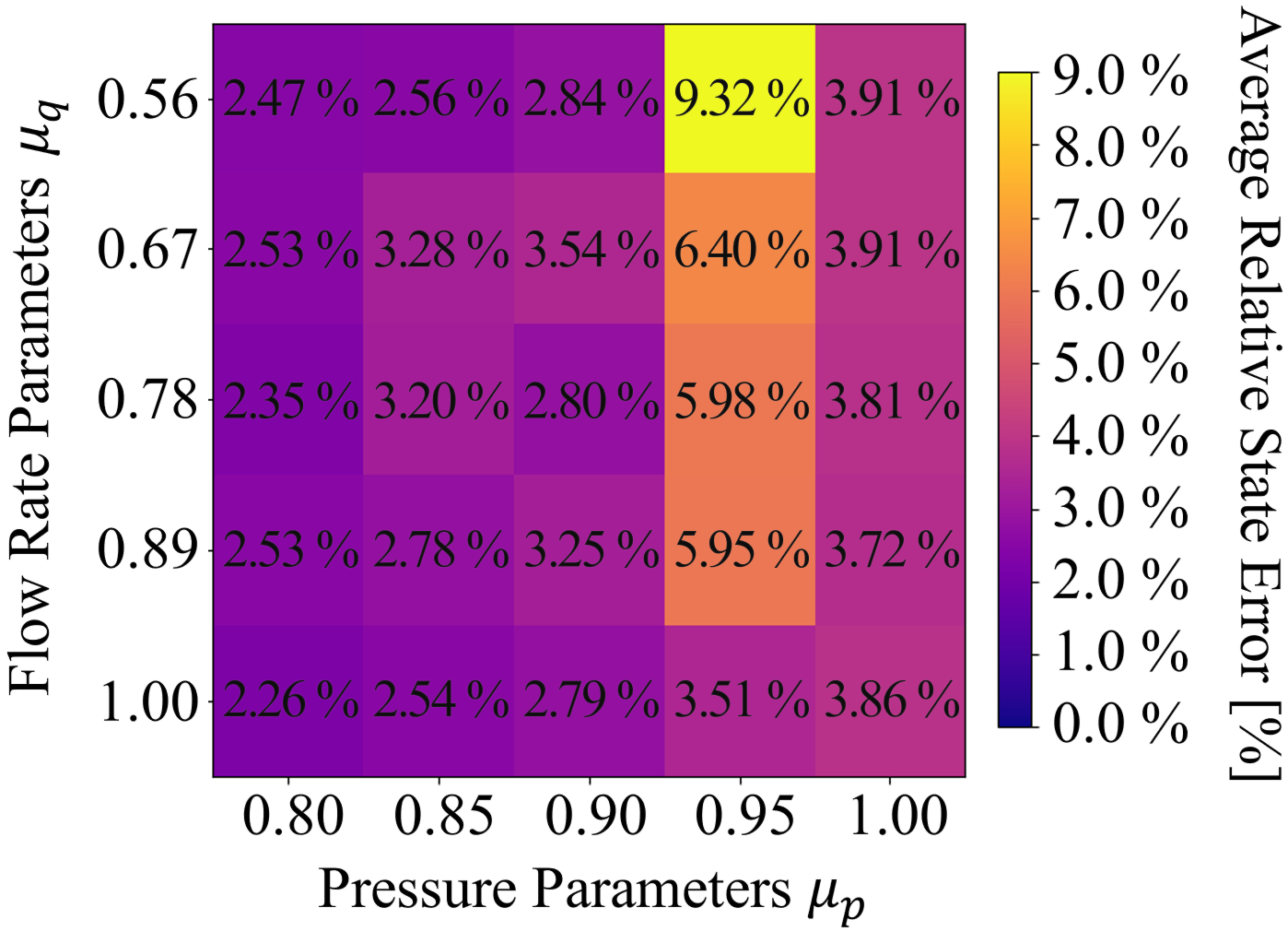}%
\label{fig:state_error_all}}
\caption{(a) Datasets used for training and testing. Of the 25 datasets, nine (36\%) are used for training, while 16 (64\%) are reserved for testing. (b) Projection error for all 25 simulations. (c) Average relative state error between the FOM and ROM predictions.}
\label{fig:errors_all_parameters}
\end{figure*}

Fig.~\ref{fig:dataset_selection} shows the selection of datasets and their corresponding parameter values used for training and testing. As discussed in Sec.~\ref{ss:FOM_details}, nine datasets, i.e., 36\% of the entire data, are used for learning the ROM.
The nine training cases are distributed across the full range of the two-dimensional parameter domain $\mu_p \in [0.80, 1.00]$ and $\mu_q \in [0.56, 1.00]$, whose convex hull forms a rectangle covering all 16 testing cases. This ensures that the surrogate model is evaluated strictly in the interpolation regime.
Fig.~\ref{fig:projection_error_all} shows the projection errors for all 25 ROMs when the global POD basis~$\bV_r$ with $r=5$ modes is used, based on nine selected training datasets. The projection error for each dataset is computed using~\eqref{eq:energy_criterion}.
The results demonstrate that our basis $\bV_r$ is well suited for both the training and testing data, with a maximum projection error of $0.74\%$ for the parameters ($\mu_q, \mu_p$) = ($0.56, 0.80$). 
Fig.~\ref{fig:state_error_all} shows the average relative state error for all 25 datasets computed using~\eqref{eq:state_error}. Among the nine training datasets, the parameter combination~($\mu_q, \mu_p$) = ($0.56, 1.00$) shows the highest average state error of $3.91$\%. For the 16 testing datasets, the lowest average error of $2.53$\% occurs at~($\mu_q, \mu_p$) = ($0.67, 0.80$) and ($0.89, 0.80$), while the highest error of~$9.32$\% is observed at~($\mu_q, \mu_p$) = ($0.56, 0.95$). Overall, the parametric OpInf ROM maintains good predictive accuracy, with errors remaining below 10\% across all 25 parameter conditions. 
\par

\begin{figure*}[!t]
\centering
\subfloat[]{\includegraphics[height=3.5in]{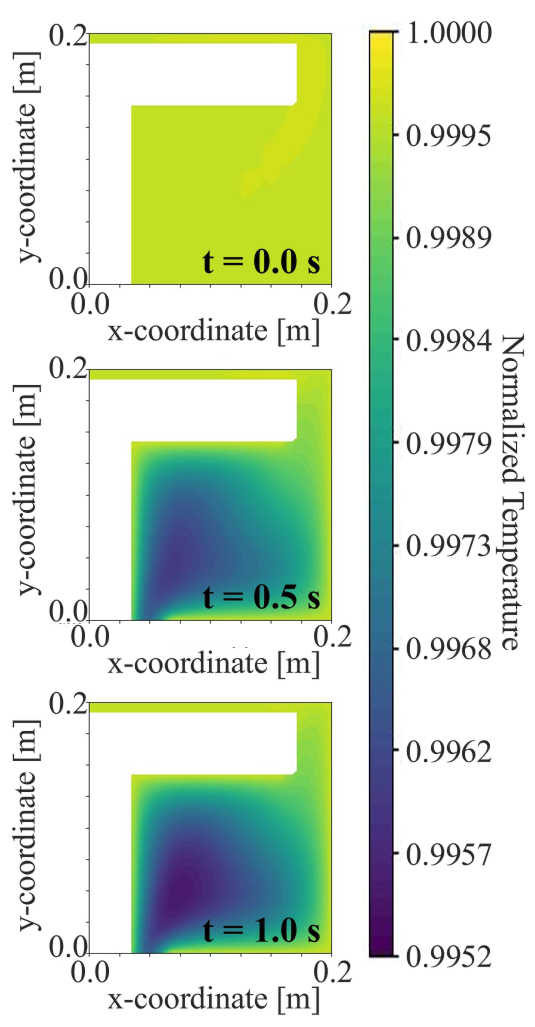}%
\label{fig:contour_FOM_temp}}
\hspace{0.5em}
\subfloat[]{\includegraphics[height=3.5in]{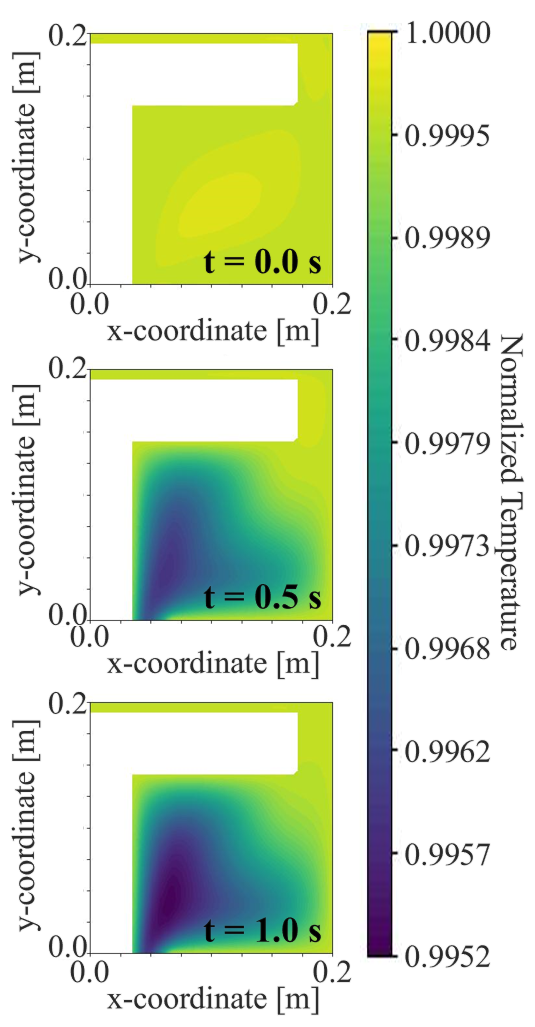}%
\label{fig:contour_ROM_temp}}
\hspace{0.5em}
\subfloat[]{\includegraphics[height=3.5in]{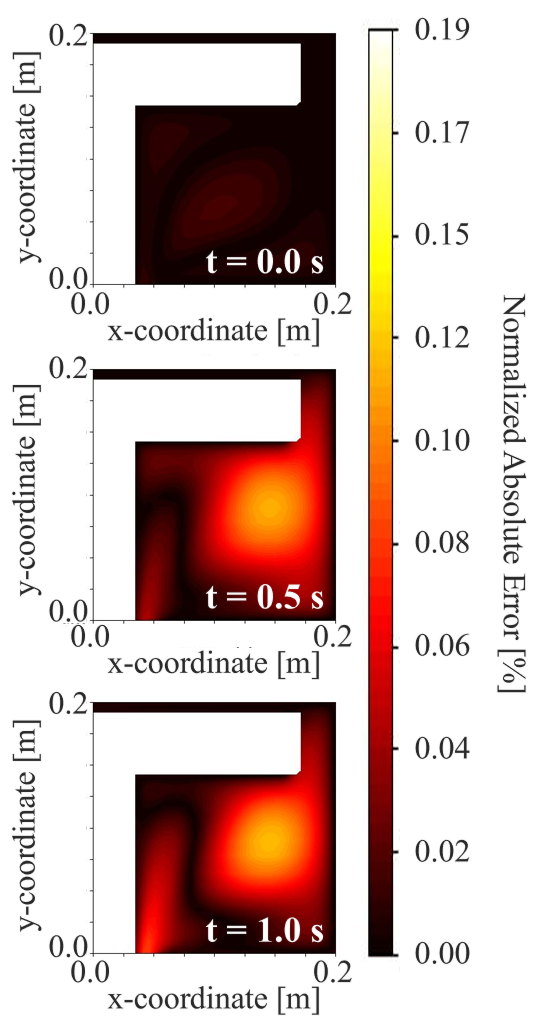}%
\label{fig:contour_ERROR_temp}}
\caption{Contour plots of the temperature in the vertical cross-section for the case with the largest error, $(\mu_p, \mu_q)=(0.56, 0.95)$. Each column represents $t=0s$, $0.5s$, and $1s$. (a) FOM, (b) ROM predictions, (c) Pointwise normalized absolute error.}
\label{fig:contour_temperature}
\end{figure*}

\begin{figure*}[!t]
\centering
\subfloat[]{\includegraphics[height=3.5in]{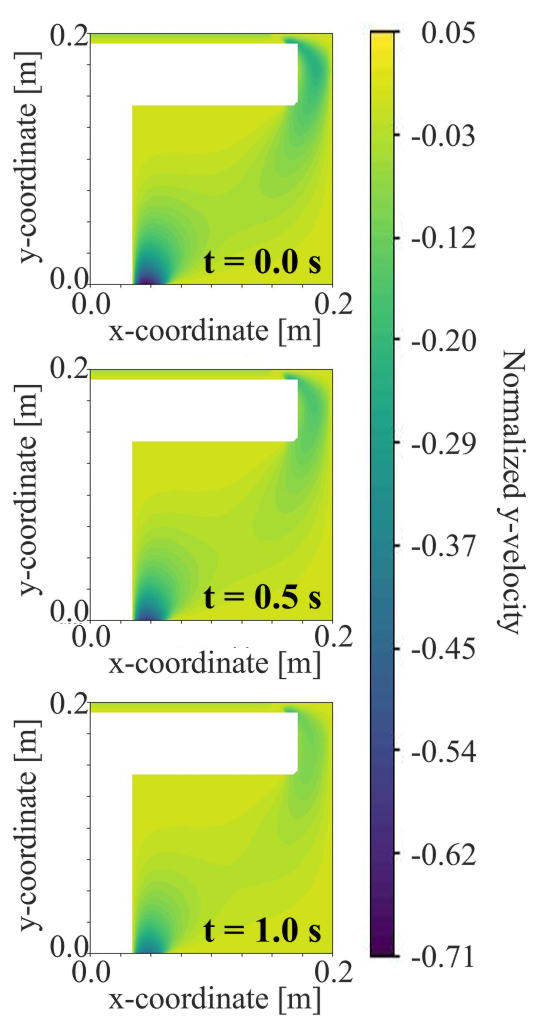}%
\label{fig:contour_FOM_y_vel}}
\hspace{0.5em}
\subfloat[]{\includegraphics[height=3.5in]{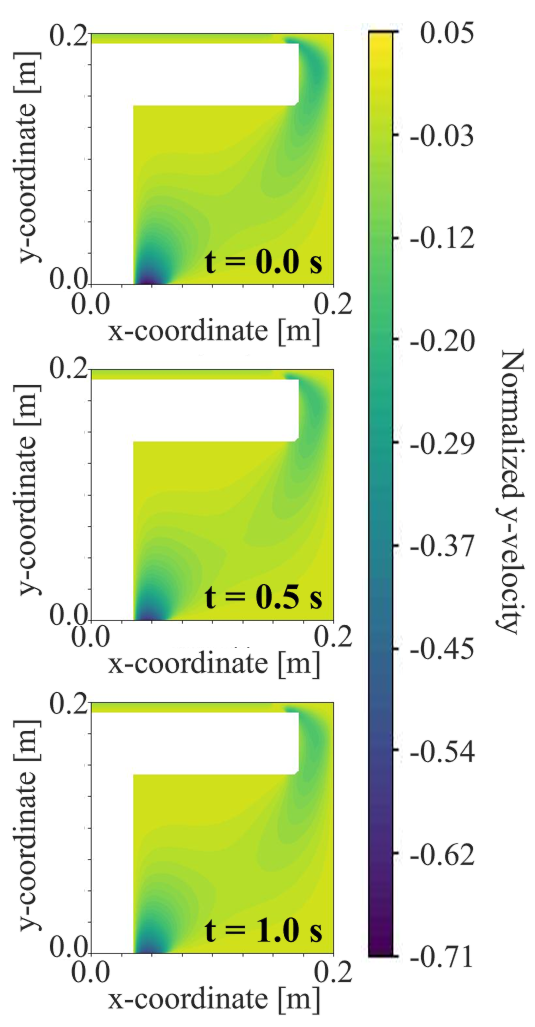}%
\label{fig:contour_ROM_y_vel}}
\hspace{0.5em}
\subfloat[]{\includegraphics[height=3.5in]{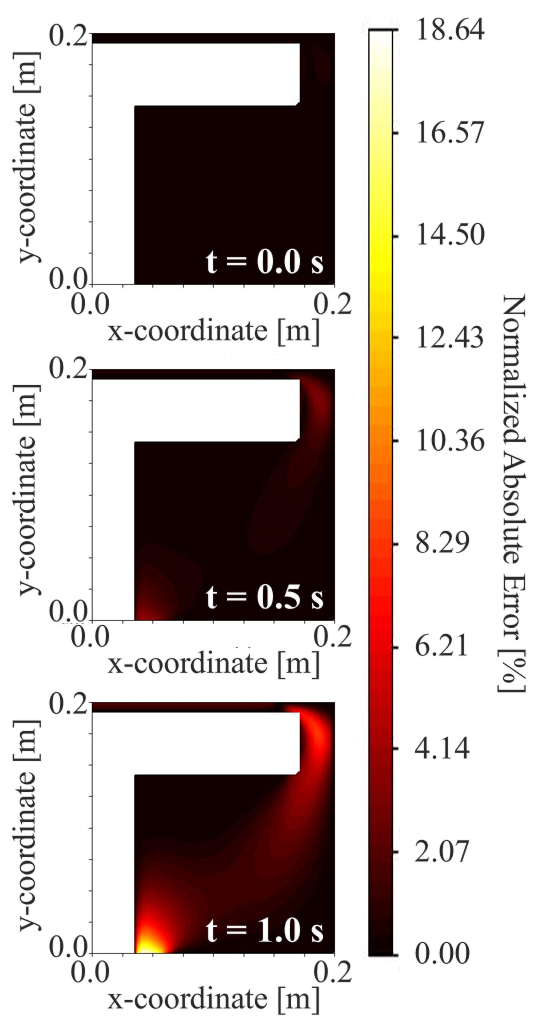}%
\label{fig:contour_ERROR_y_vel}}
\caption{Contour plots of the y-velocity $v_y$ in the vertical cross-section for the case with the largest error, $(\mu_p, \mu_q)=(0.56, 0.95)$. Each column represents $t=0s$, $0.5s$, and $1s$. (a) FOM, (b) ROM predictions, (c) Pointwise normalized absolute error.}
\label{fig:contour_y_velocity}
\end{figure*}

\begin{figure*}[!t]
\centering
\subfloat[]{\includegraphics[height=3.5in]{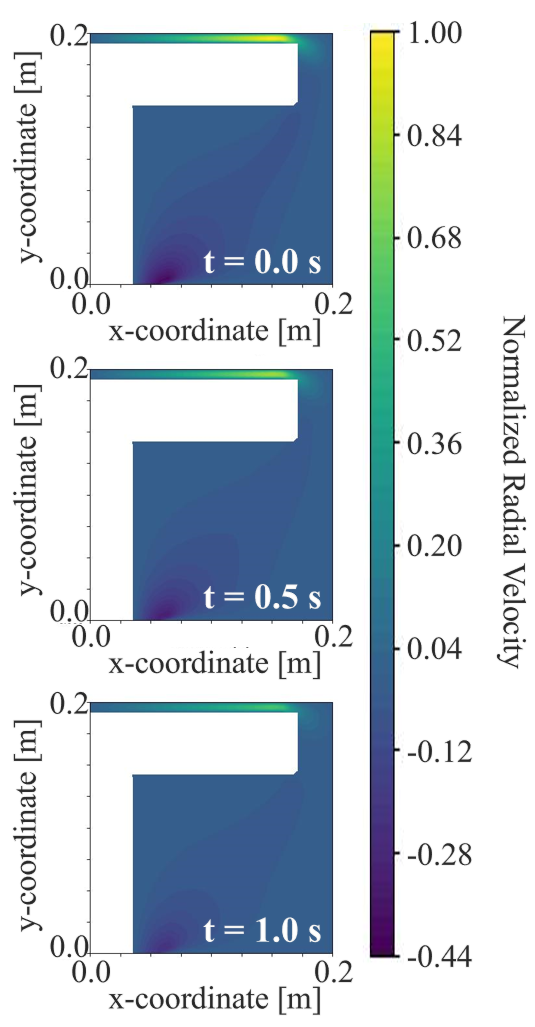}%
\label{fig:contour_FOM_radial_vel}}
\hspace{0.5em}
\subfloat[]{\includegraphics[height=3.5in]{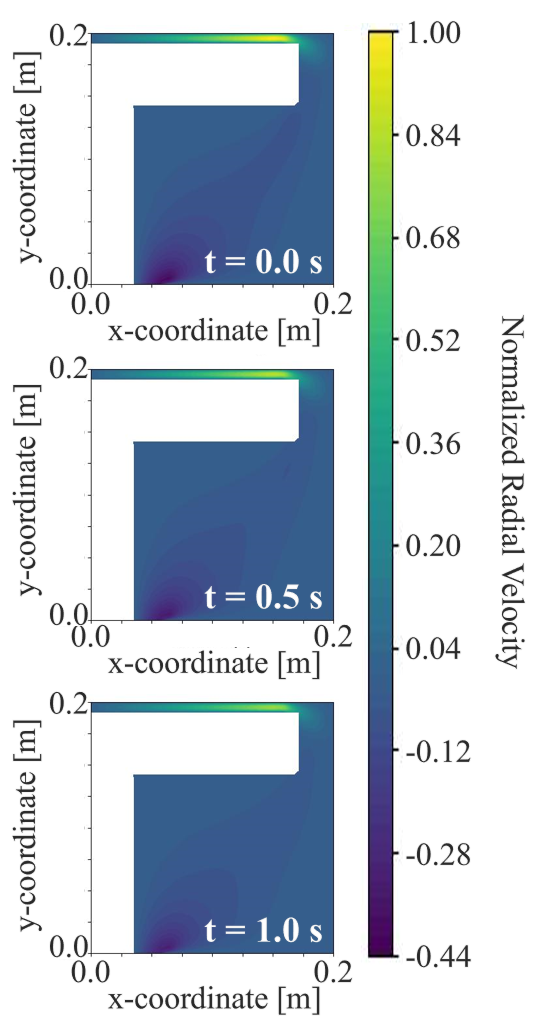}%
\label{fig:contour_ROM_radial_vel}}
\hspace{0.5em}
\subfloat[]{\includegraphics[height=3.5in]{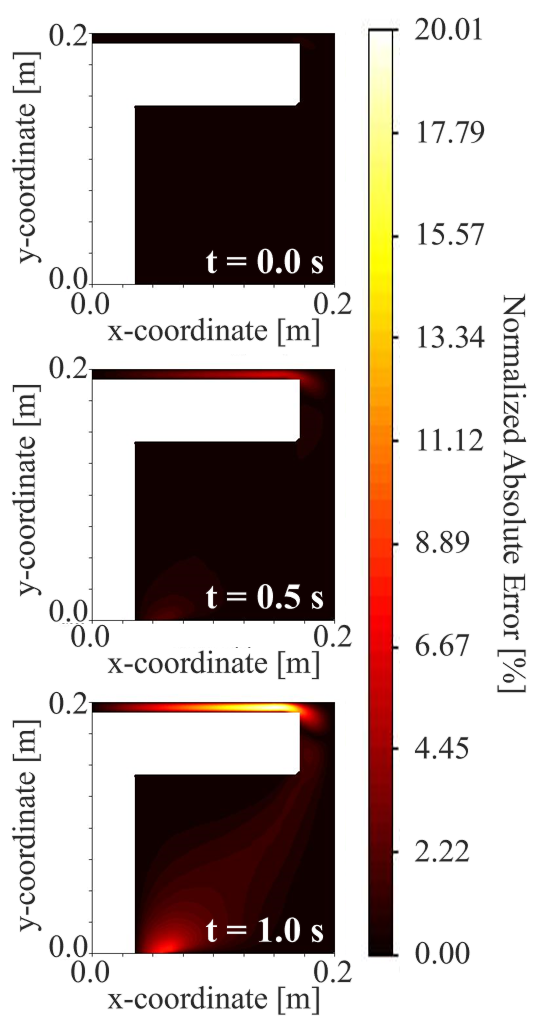}%
\label{fig:contour_ERROR_radial_vel}}
\caption{Contour plots of the radial velocity $v_r=\sqrt{v_x^2 + v_z^2}$ in the vertical cross-section for the case with the largest error, $(\mu_p, \mu_q)=(0.56, 0.95)$. Each column represents $t=0s$, $0.5s$, and $1s$. (a) FOM, (b) ROM predictions, (c) Pointwise normalized absolute error.}
\label{fig:contour_radial_velocity}
\end{figure*}

Figs.~\ref{fig:contour_temperature},~\ref{fig:contour_y_velocity}, and~\ref{fig:contour_radial_velocity} show the flow fields of both the FOM and ROM predictions over time, along with their corresponding pointwise errors, defined as
\begin{equation*}
    \Upsilon_{j,k}^i = \frac{\left| \bS_{j,k}^{\text{FOM},i} - \bS_{j,k}^{\text{ROM},i} \right|}{\max \left( \left|\bS_{j,k}^{\text{FOM},i}\right| \right)},
\end{equation*}
where $j$ and $k$ are the spatial and temporal indices, respectively. Thus, the error $\Upsilon_{j,k}^i$ is the normalized absolute error of each unscaled variable over the entire spatiotemporal domain.
All figures represent a vertical cross-section of the $xy$-vertical plane where $z=0$~mm, as shown in Fig.~\ref{fig:computational_domain}. To effectively compare the differences in scale, the velocity components are normalized by the maximum absolute value between $v_y$ and $v_r$. The case presented corresponds to the parameters~$(\mu_p, \mu_q)=(0.56, 0.95)$, which yields the largest average relative state error (see Fig.~\ref{fig:state_error_all}).   
We omit the pressure results since they are nearly uniform and exhibit no significant gradients across the entire vertical cross-sectional area. However, we note that the ROM captures the FOM pressure behavior, with a maximum pointwise error of 0.429\% at $t=1s$. 
Fig.~\ref{fig:contour_temperature} shows that the ROM also predicts the temperature with high accuracy, producing a maximum error of $0.19\%$. Although this maximum error occurs specifically under the wafer heater, the ROM captures the temperature distribution across the entire spatiotemporal domain, especially given the small variation in its scale.  \par

\begin{figure}[!t]
\centering
\includegraphics[height=2.7in]{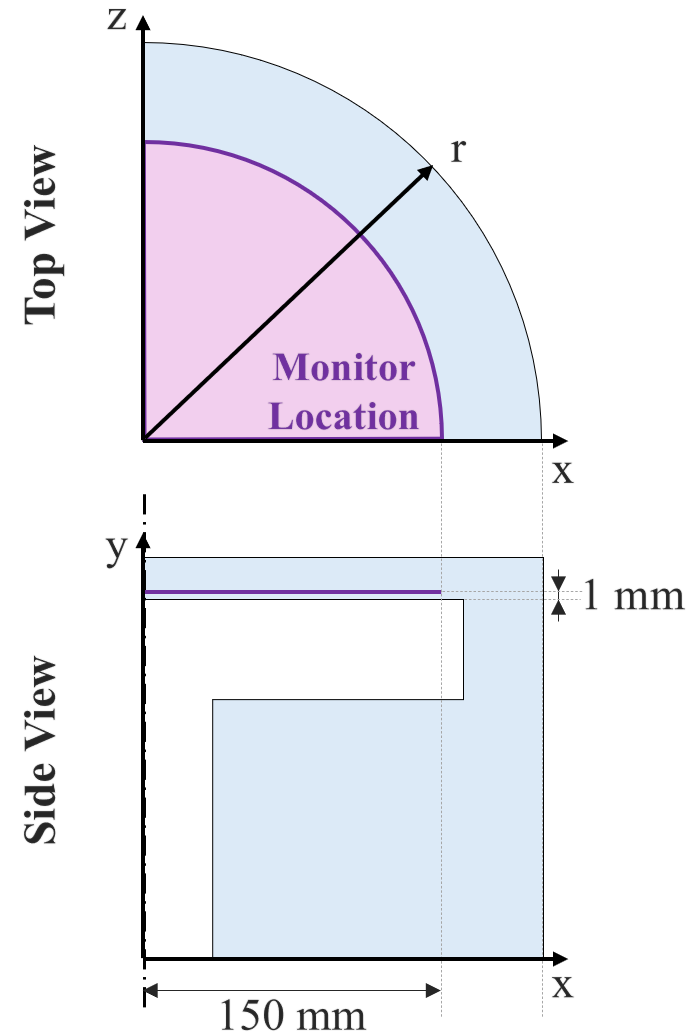}
\caption{Monitor location: a quarter circular $xz$-plane with a 150mm radius, positioned 1~mm above the wafer surface.}
\label{fig:monitor_location}
\end{figure}

The velocity components, however, exhibit relatively high errors, as shown in Fig.~\ref{fig:contour_y_velocity} and~\ref{fig:contour_radial_velocity}. In particular, the y-directional and radial velocities show maximum pointwise errors of~$18.64\%$ and~$20.01\%$, respectively.
This case corresponds to the parameter combination~$(\mu_p, \mu_q)=(0.56, 0.95)$, which produces the highest average error of~$9.32\%$ across pressure, temperature, y-velocity, and radial velocity, indicating that the velocity errors are the dominant contributors to this overall error.
To analyze more in detail, the comparison between Figs.~\ref{fig:contour_FOM_y_vel} and~\ref{fig:contour_ROM_y_vel} shows that the ROM accurately captures the overall direction of the y-velocity of the purging flow. However, the error contours in Fig.~\ref{fig:contour_ERROR_y_vel} highlight the difference in magnitude between the FOM and ROM predictions, with a maximum error of 18.64\% near the chamber outlet at the final time step ($t=1s$).
Similarly, Figs.~\ref{fig:contour_FOM_radial_vel} and~\ref{fig:contour_ROM_radial_vel} demonstrate that the ROM captures the overall pattern of the radial velocity field. However, Fig.~\ref{fig:contour_ERROR_radial_vel} shows that the errors increase toward the edge of the wafer surface at $t=1$s.  
From a practical perspective, accurate predictions of flow fields in the specific region that directly impacts particle contamination are especially important. During the PECVD process, after the plasma is deactivated, particles are carried by the purging flow. Since the level of contamination is determined by the amount of particles deposited on the wafer surface, the flow dynamics directly above the wafer surface are critical for particle contamination control. Thus, the monitor plane is defined as the $xz$-plane at the height of 1\text{mm} above the wafer surface (see~Fig.~\ref{fig:monitor_location}), as this region is most susceptible to particle contamination and therefore most critical for prediction accuracy. \par

\begin{figure}[!ht]
    \centering
    \begin{subfigure}[b]{0.3\linewidth}
        \centering
        \includegraphics[height=8.7cm, keepaspectratio]{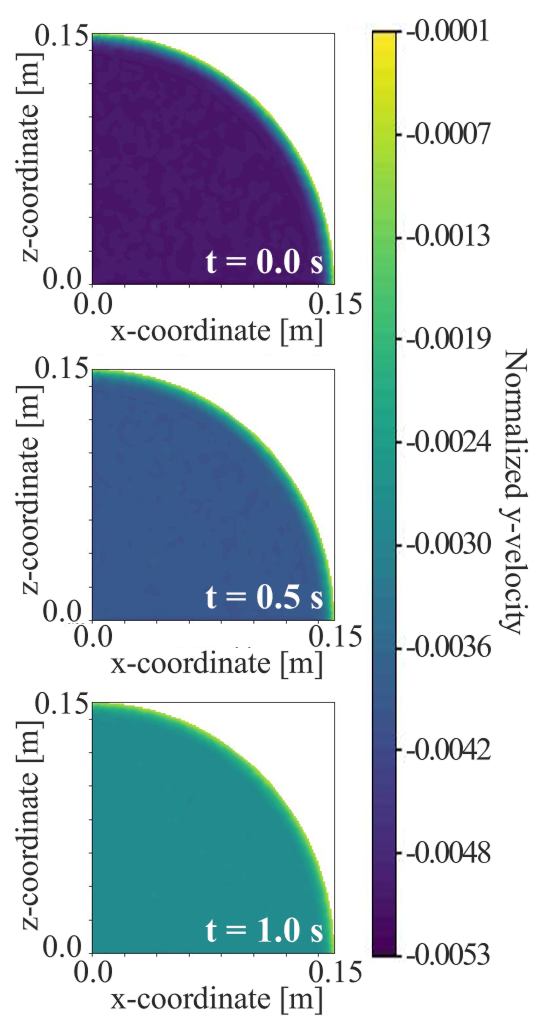}
        \caption{}
        \label{fig:mon_contour_FOM_y_vel}
    \end{subfigure}
    \hspace{0cm}
    \begin{subfigure}[b]{0.3\linewidth}
        \centering
        \includegraphics[height=8.7cm, keepaspectratio]{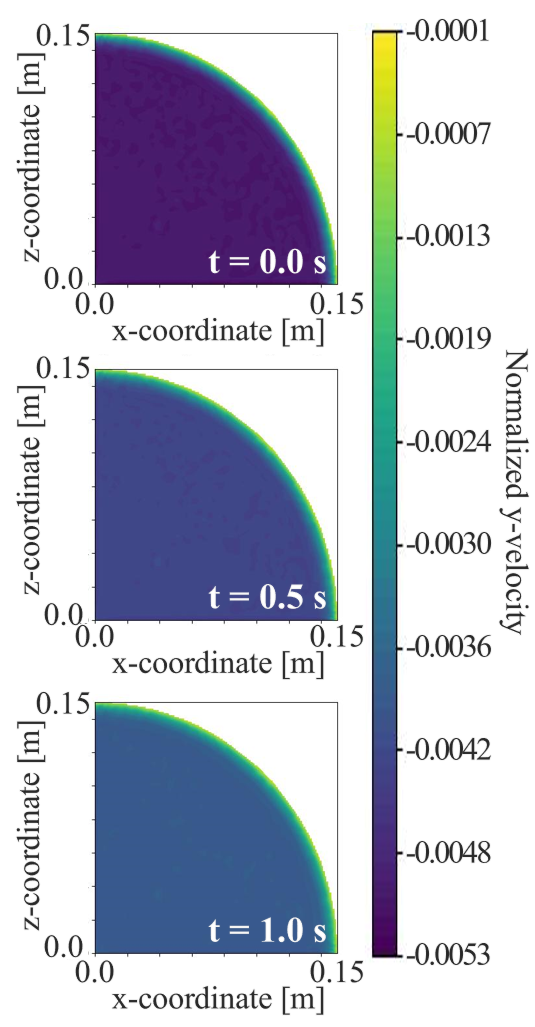}
        \caption{}
        \label{fig:mon_contour_ROM_y_vel}
    \end{subfigure}        
    \hspace{0cm}
    \begin{subfigure}[b]{0.3\linewidth}
        \centering
        \includegraphics[height=8.7cm, keepaspectratio]{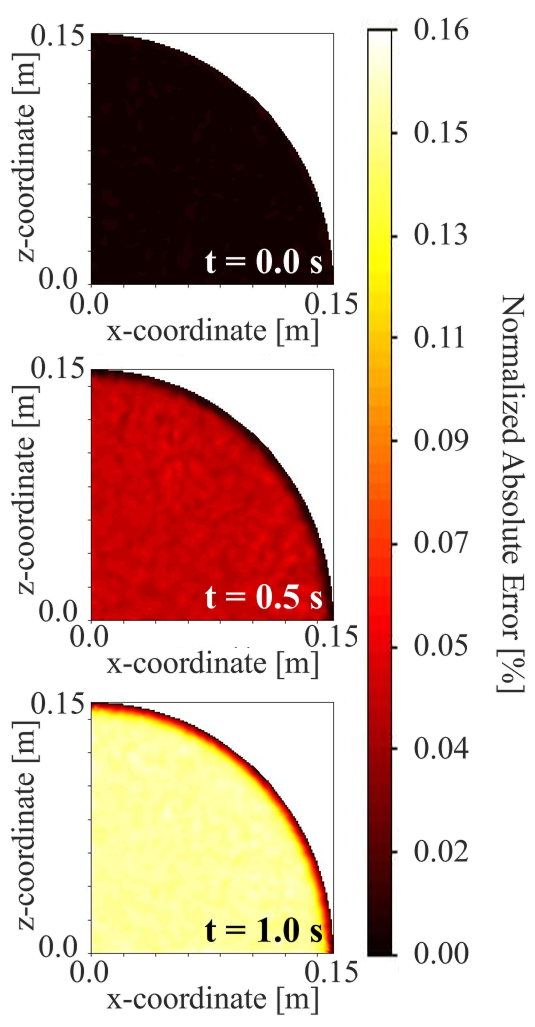}
        \caption{}
        \label{fig:mon_contour_ERROR_y_vel}
    \end{subfigure}    
    \vspace{-0.5em}
    \caption{Y-velocity contour plots at the monitor location for the case with the largest error, $(\mu_p, \mu_q)=(0.56, 0.95)$. Each column represents $t=0s$, $0.5s$, and $1s$. (a) FOM, (b) ROM predictions, (c) Pointwise normalized absolute error.
    }
    \label{fig:mon_contour_y_velocity}
\end{figure}

\begin{figure}[!ht]
    \centering
    \begin{subfigure}[b]{0.3\linewidth}
        \centering
        \includegraphics[height=8.7cm, keepaspectratio]{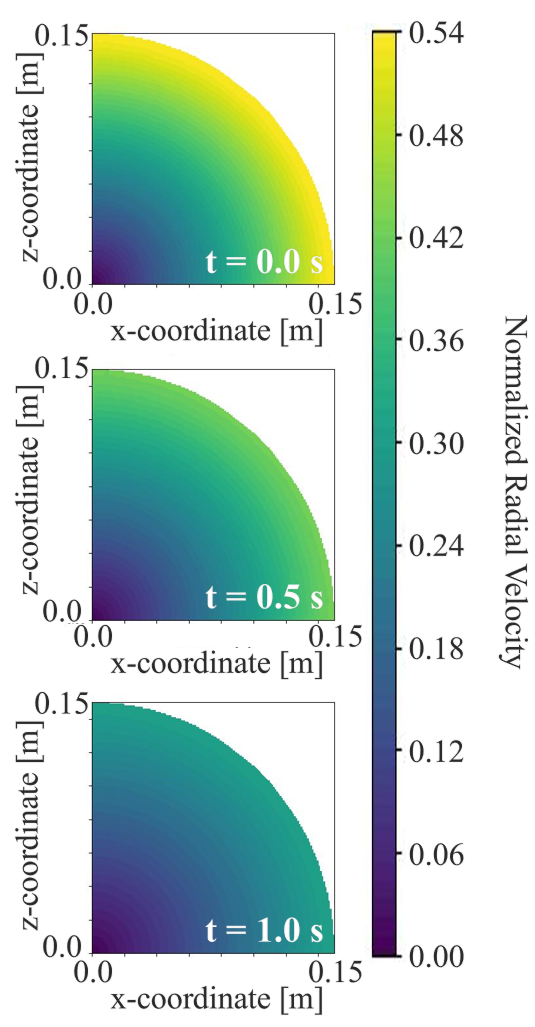}
        \caption{}
        \label{fig:mon_contour_FOM_radial_vel}
    \end{subfigure}
    \hspace{0cm}
    \begin{subfigure}[b]{0.3\linewidth}
        \centering
        \includegraphics[height=8.7cm, keepaspectratio]{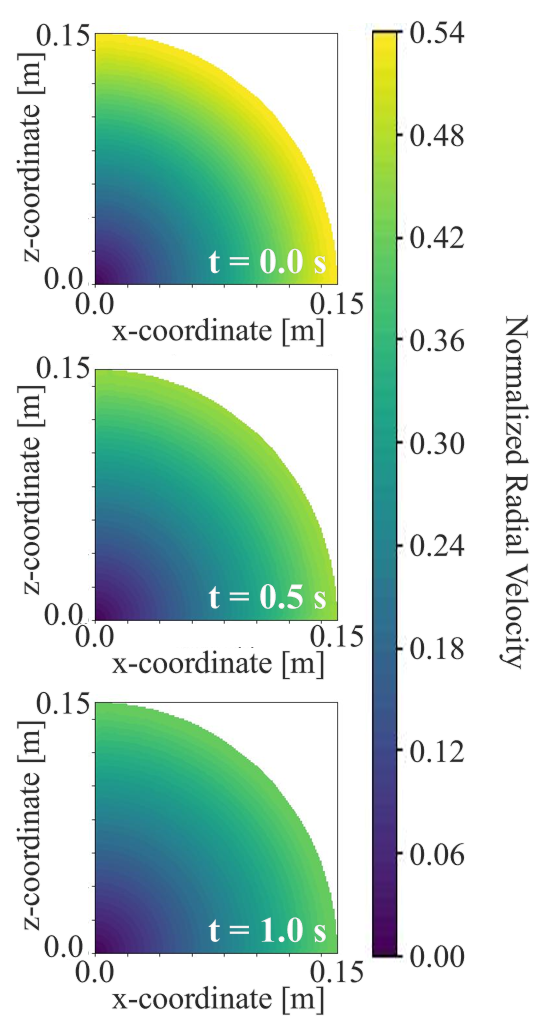}
        \caption{}
        \label{fig:mon_contour_ROM_radial_vel}
    \end{subfigure}        
    \hspace{0cm}
    \begin{subfigure}[b]{0.3\linewidth}
        \centering
        \includegraphics[height=8.7cm, keepaspectratio]{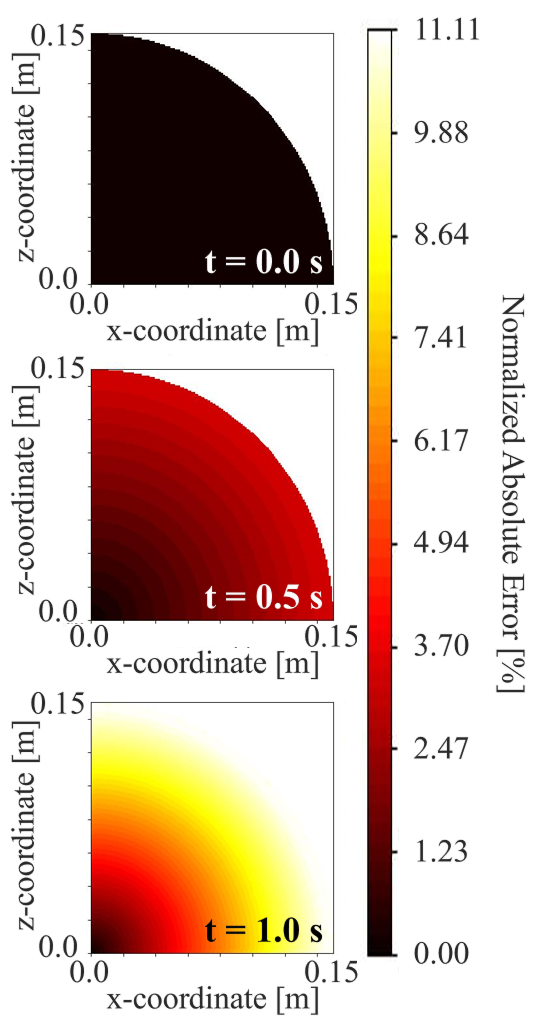}
        \caption{}
        \label{fig:mon_contour_ERROR_radial_vel}
    \end{subfigure}    
    \vspace{-0.5em}
    \caption{
    Radial velocity contour plots at the monitor location for the case with the largest error, $(\mu_p, \mu_q)=(0.56, 0.95)$. Each column represents $t=0s$, $0.5s$, and $1s$. (a) FOM, (b) ROM predictions, (c) Pointwise normalized absolute error.
    }
    \label{fig:mon_contour_all}
\end{figure}

We compare the y-velocity and radial velocity flow fields at the monitor location for the case with the largest error, $(\mu_p, \mu_q)=(0.56, 0.95)$, in Fig.~\ref{fig:mon_contour_all}. The velocities are normalized by the maximum absolute value between $v_y$ and $v_r$, highlighting that the radial velocity has a much larger magnitude than the y-velocity.
As shown in Fig.~\ref{fig:mon_contour_FOM_y_vel}, the y-velocity is nearly uniform across the monitor surface, except at the edge, where it changes due to the acceleration of the purging flow. As shown in Fig.~\ref{fig:mon_contour_ROM_y_vel}, the ROM captures this behavior well, which is also partially observed in the area on top of the heater in Fig.~\ref{fig:contour_ERROR_y_vel}. Moreover, Fig.~\ref{fig:mon_contour_ERROR_y_vel} demonstrates that the ROM provides accurate y-velocity predictions across the entire surface, with a maximum error of 0.16\%.
In contrast, the radial velocity exhibits higher errors, as shown in Fig.~\ref{fig:mon_contour_ERROR_radial_vel}, especially near the wafer edge at $t=1$s.
Taken together, the ROM captures the overall trend of radial velocity fields on the wafer surface well with a small number of POD modes~($r=5$), but it falls short in accurately predicting very localized features where the FOM exhibits sharp gradients near the edge. This is a common feature of global linear modal approximations. 
While incorporating---heuristically---more local mode shapes or nonlinear manifolds could potentially alleviate this issue, such extensions significantly increase the complexity of the OpInf problem \cite{GEELEN2023115717}. For further discussions on both intrusive and nonintrusive localized ROM frameworks, see \cite{amsallem2012nonlinear, geelen2022localized, buhr2020localized}.
\par
From a purging perspective, the two dominant velocity components of the purging flow serve distinct roles: the y-velocity vertically directs particles toward the wafer surface, while the radial velocity transports them laterally from the center toward the edge and ultimately to the exhaust outlet. Consequently, particles near the center must travel a longer distance to reach the edge, driven by the combined action of the radial and y-velocity, making this region inherently more susceptible to contamination.
Near the wafer center, where $(x, z) = (0,0)[\textnormal{m}]$, the y-velocity reaches its highest magnitude while the radial velocity is at its minimum, yielding a low radial-to-y velocity ratio (see Fig.~\ref{fig:mon_contour_all}). This flow configuration promotes particle deposition onto the wafer surface, increasing contamination risk. Conversely, near the wafer edge, the radial velocity is dominant, producing a high radial-to-y velocity ratio that efficiently sweeps particles toward the exhaust outlet, reducing local contamination risk.
The spatial ROM prediction errors are consistent with these contamination risk regions. 
Despite the maximum y-velocity error of 18.64\% occuring near the outlet (bottom of the chamber), the error within the wafer center region remains below 0.16\%, see Figs.~\ref{fig:contour_ERROR_y_vel} and~\ref{fig:mon_contour_ERROR_y_vel}. This confirms that the elevated y-velocity errors are spatially confined to the outlet vicinity, a region of negligible contamination risk.
Similarly, while the radial velocity exhibits its highest error of 11.11\% near the wafer edge, this region poses relatively lower contamination risk due to the dominant radial flow that promotes efficient particle purging. 
Moreover, as shown in~Fig.~\ref{fig:mon_contour_ERROR_radial_vel}, the radial velocity error decreases progressively toward the wafer center, the region of the highest contamination risk.
This further confirms that the ROM maintains high predictive accuracy precisely in the area that is most critical to particle contamination.

\subsection{Computational Speedup of the OpInf ROM vs.~CFD} \label{ss:Computing_efficiency}
We measure the CPU time of the OpInf ROM to compare it with that of the FOM to evaluate computational efficiency.
The OpInf ROM simulation consists of two stages: the offline stage, where the model is constructed, and the online stage, where predictions are made using the learned ROM.
To assess the computational efficiency of each stage and identify potential computational bottlenecks, we measure the CPU time separately for the offline and online stages. 
To account for variations in ROM simulation time, we compute the average online cost over 20 simulations for all 25 ROMs. Both FOM and ROM computations are done using 20 CPU cores on a 13th-generation Intel Core i5-13600K processor (3.50 GHz).
\Cref{tab:computing_time} compares the measured CPU time and shows the computational speedup of the ROM. To account for variations in ROM simulation time depending on various factors, we take the average of the online costs over 20 ROM simulations for 25 parameter combinations. The ROM construction in the offline stage takes 1491.27 minutes. Note that the model is constructed with nine parameter combinations (see~ Fig.~\ref{fig:dataset_selection}). However, the online stage, which predicts the flow for all 25 parameter combinations (including both training and testing cases), takes only 1.24 minutes on average over 20 simulations. This corresponds to about a 142-times speedup compared to the FOM computation time of 176.88 minutes. In other words, for many-query situations where we need to evaluate the ROM more than $1491.27/176.88 \approx 8$ times, it pays off to construct a ROM instead of running the FOM. 
\begin{table}[!t]    
\caption{CPU time comparison. The online speedup factor is computed by dividing the FOM simulation cost by the online ROM simulation cost.}  \label{tab:computing_time}
\centering
\begin{tabular}{c|c c c|c}
& \textbf{FOM} & \textbf{ROM} & \textbf{ROM} & \textbf{Online} \\ 
&  & \textbf{offline} & \textbf{online} & \textbf{speedup factor} \\
\hline
\textbf{CPU time [min]} & 176.88  & $1491.27$ & $1.24$ & $142.65$
\end{tabular}
\end{table}

\section{Conclusion} \label{sec:conclusion}
We simulated the purging process in the PECVD chamber using OpInf ROMs, which leverage a data-driven approach to model dynamical systems with complex multiscale and nonlinear behavior.
The OpInf framework interpolated between nine ROMs constructed using different argon flow rates at the inlet and outlet pressures in the purging process. This enabled the prediction of the system's behavior for 16 unseen parameter combinations.
The OpInf ROMs maintained high predictive accuracy for the pressure and temperature across the entire spatiotemporal domain. 
While velocity predictions showed relatively higher errors near the wafer edge, the ROM captured critical flow dynamics in the wafer center region, which is sensitive to particle contamination.
In general, the parametric OpInf ROMs showed good predictive accuracy across 25 parameter combinations, with a maximum error of 9.32\%.
This approach substantially lowered computational costs, achieving a 142-fold speedup in online computations compared to the FOM CFD simulation.
Overall, this fast and accurate predictive model for the flow field in the PECVD purging process enables the prediction of particle movement within the PECVD chamber, making it a useful tool for particle contamination control. 
This capability will be particularly crucial in the fast-paced and precision-driven semiconductor manufacturing environment.
\par
There are several interesting directions for future work. First, the performance of the proposed approach may vary depending on the choice of training and testing datasets. We selected a low training-to-testing ratio (9 train/16 test) to demonstrate the predictive accuracy of the ROM, see Fig.~\ref{fig:dataset_selection}. While adopting a more conventional split (e.g., 20 training and five testing datasets) would likely improve overall accuracy, our goal was to highlight the model’s interpolation performance under limited training. An interesting direction for future work would be to explore the ROM’s extrapolation capability beyond the training domain---for example, in cases where $\mu_p < 0.80$ or $\mu_q < 0.56$.
Second, we note that the contamination-related analysis presented in this work is inferential, relying on flow-field predictions, rather than explicit particle tracking or contamination metrics. Specifically, the surrogate model is validated through the accuracy of the velocity components most directly governing particle transport---the radial and y-velocity fields in the contamination-critical wafer center region---rather than through direct simulation of particle trajectories. A natural extension would be to couple particle-tracking models with the flow-field predicting ROMs. This integration would enable direct simulation of particle transport and deposition, and comparison with particle measurement data from the PECVD chamber, thereby offering a more comprehensive framework for contamination risk prediction.
Lastly, it would also be interesting to investigate whether this approach would work for faster transient phenomena, such as fast atomic layer deposition steps, where the flow dynamics vary significantly on much smaller time scales.
    
\section*{Acknowledgments}
This research was financially supported by Samsung Electronics Co., Ltd., under award 30312263 for the project ``Reduced-Order Modeling for Real-Time Simulation of Flow Phenomena in Semiconductor Manufacturing" and gift funds from ASML US, LP.

\bibliography{references}
\bibliographystyle{abbrv}

\end{document}